\documentclass[reqno,11pt]{amsart}
\usepackage{amssymb,enumerate}
\usepackage{amsmath}
\usepackage{mathrsfs}
\usepackage{dsfont} %pour avoir indicatrice avec \mathds{1}
\usepackage{hyperref}
\usepackage{dsfont}
\usepackage[usenames]{color}

\allowdisplaybreaks

\begin{document}
%\renewcommand{\footnote}{}
%\pagenumbering{roman}
%\thispagestyle{empty}%\quad\newpage
%\thispagestyle{empty}
%\nonstopmode
%**************************************************************************

%%%%%%%%%%%%%%%%%%%%%%%%%%%%%%%%%%%%%%%%%%%%%%%%%%%%%%%%%%%%%%%%%%%%%%%%%%%%%%%%%
\renewcommand{\AA}{\mathcal{A}}
\renewcommand{\div}{{\rm div \,}}
\newcommand{\BB}{\mathcal{B}}
\newcommand{\CC}{\mathcal{C}}
\newcommand{\DD}{\mathcal{D}}
\newcommand{\EE}{\mathcal{E}}
\newcommand{\FF}{\mathcal{F}}
\newcommand{\GG}{\mathcal{G}}
\newcommand{\HH}{\mathcal{H}}
\newcommand{\II}{ I\hspace{-0.5mm}I}
\newcommand{\JJ}{\mathcal{J}}
\newcommand{\KK}{\mathcal{K}}
\newcommand{\LL}{\mathcal{L}}
\newcommand{\MM}{\mathcal{M}}
\newcommand{\NN}{\mathcal{N}}
\newcommand{\OO}{\mathcal{O}}
\newcommand{\PP}{\mathcal{P}}
\newcommand{\QQ}{\mathcal{Q}}
\newcommand{\RR}{\mathcal{R}}
\renewcommand{\SS}{\mathcal{S}}
\newcommand{\TT}{\mathcal{T}}
\newcommand{\UU}{\mathcal{U}}
\newcommand{\VV}{\mathcal{V}}
\newcommand{\WW}{\mathcal{W}}
\newcommand{\XX}{\mathcal{X}}
\newcommand{\YY}{\mathcal{Y}}
\newcommand{\ZZ}{\mathcal{Z}}

%Mengensystem

\newcommand{\ds}{\displaystyle}

\newcommand{\A}{\mathbb{A}}
\newcommand{\B}{\mathbb{B}}
\newcommand{\C}{\mathbb{C}}
\newcommand{\D}{\mathbb{D}}
\newcommand{\E}{\mathbb{E}}
\newcommand{\F}{\mathbb{F}}
\newcommand{\G}{\mathbb{G}}
\renewcommand{\H}{\mathbb{H}}
\newcommand{\I}{\mathbb{I}}
\newcommand{\J}{\mathbb{J}}
\newcommand{\K}{\mathbb{K}}
\renewcommand{\L}{\mathbb{L}}
\newcommand{\M}{\mathbb{M}}
\newcommand{\N}{\mathbb{N}}
\renewcommand{\O}{\mathbb{O}}
\renewcommand{\P}{\mathbb{P}}
\newcommand{\Q}{\mathbb{Q}}
\newcommand{\R}{\mathbb{R}}
\renewcommand{\S}{\mathbb{S}}
\newcommand{\T}{\mathbb{T}}
\newcommand{\U}{\mathbb{U}}
\newcommand{\V}{\mathbb{V}}
\newcommand{\W}{\mathbb{W}}
\newcommand{\X}{\mathbb{X}}
\newcommand{\Y}{\mathbb{Y}}
\newcommand{\Z}{\mathbb{Z}}

%griechisch

\newcommand{\al}{\alpha}
\newcommand{\be}{\beta}
\newcommand{\ga}{\gamma}
\newcommand{\de}{\delta}
\newcommand{\ep}{\varepsilon}
\newcommand{\ze}{\zeta}
\newcommand{\et}{\eta}
\newcommand{\vth}{\vartheta}
\renewcommand{\th}{\theta}
\newcommand{\io}{\iota}
\newcommand{\ka}{\kappa}
\newcommand{\la}{\lambda}
\newcommand{\rh}{\rho}
\newcommand{\si}{\sigma}
\newcommand{\ta}{\tau}
\newcommand{\up}{\upsilon}
\newcommand{\ph}{\varphi}
\newcommand{\ch}{\chi}
\newcommand{\ps}{\psi}
\newcommand{\om}{\omega}

\newcommand{\Ga}{\Gamma}
\newcommand{\De}{\Delta}
\newcommand{\Th}{\Theta}
\newcommand{\La}{\Lambda}
\newcommand{\Si}{\Sigma}
\newcommand{\Up}{\Upsilon}
\newcommand{\Ph}{\xi}
\newcommand{\Om}{\Omega}

\newcommand{\BR}{\color{red}}
\newcommand{\BL}{\color{blue}}
\newcommand{\ER}{\color{black}}

%mathematische Symbole
%%Pfeile
\newcommand{\inj}{\hookrightarrow}
\newcommand{\stetein}{\overset{s}{\hookrightarrow}}
\newcommand{\dichein}{\overset{d}{\hookrightarrow}}
\newcommand{\pa}{\partial}
\newcommand{\re}{\restriction}
\newcommand{\tief}{\downharpoonright}

%%Diverses
\newcommand{\bra}{\langle}
\newcommand{\ket}{\rangle}
\newcommand{\bs}{\backslash}
\newcommand{\divv}{\operatorname{div}}
\newcommand{\Dt}{\frac{\mathrm d}{\mathrm dt}}

%%Mengenoperationen
\newcommand{\sm}{\setminus}
\newcommand{\es}{\emptyset}

\newtheorem{theorem}{Theorem}%[section]
\newtheorem{corollary}{Corollary}
\newtheorem*{main}{Main Theorem}
\newtheorem{lemma}[theorem]{Lemma}
\newtheorem{proposition}{Proposition}
\newtheorem{conjecture}{Conjecture}
\newtheorem*{problem}{Problem}
\theoremstyle{definition}
\newtheorem{definition}[theorem]{Definition}
\newtheorem{remark}{Remark}
\newtheorem*{notation}{Notation}

%%Additional French shortcuts
\newcommand{\cqfd}{\begin{flushright}\vspace*{-3mm}$\Box $\vspace{-2mm}\end{flushright}}
\newcommand{\saut}{\vspace*{1mm}\\ \nd }
\newcommand{\on}{\mbox{ on }}
\newcommand{\with}{\mbox{ with }}
\newcommand{\nd}{\noindent}
\newcommand{\eps}{\varepsilon}
\newcommand{\limit}{\underset{n\rightarrow +\infty}{\longrightarrow}}

%**************************************************************************
\sloppy
%**************************************************************************
\title[]
{Global Wellposedness for a class\\ of Reaction-Advection-Anisotropic-Diffusion Systems}
\author[D. Bothe]{Dieter Bothe}
\address{Center of Smart Interfaces, Technische Universit\"at Darmstadt, Alarich-Weiss-Str.\ 10, 64287 Darmstadt, Germany}
\email{bothe@csi.tu-darmstadt.de}

\author[A. Fischer]{Andr\'e Fischer}
\address{Center of Smart Interfaces, Technische Universit\"at Darmstadt, Alarich-Weiss-Str.\ 10, 64287 Darmstadt, Germany}
\email{fischer@csi.tu-darmstadt.de}

\author[M. Pierre]{Michel Pierre}
\address{Ecole Normale Sup\'erieure de Rennes, IRMAR, UEB, Campus de Ker Lann, 35170 Bruz, France}
\email{michel.pierre@ens-rennes.fr}

\author[G. Rolland]{Guillaume Rolland}
\address{Ecole Normale Sup\'erieure de Rennes, IRMAR, UEB, Campus de Ker Lann, 35170 Bruz, France}
\email{guillaume.rolland@ens-cachan.org}

\date{\today}
%\subjclass[2010]{Primary: 76E20, 76U05; Secondary: 35Q30.}
% \keywords{Ekman spiral, Navier-Stokes equations, nonlinear instability.}
%\thanks{class-stefan75.tex}
\thispagestyle{empty}
%\pagestyle{myheadings}
%\markboth{\today}{\today}
%\setlength{\parindent}{0mm}
%\setlength{\parskip}{10cm}
\parskip0.5ex plus 0.5ex minus 0.5ex
%\bibliographystyle{plain}
%\bibliographystyle{alpha}
%\bibliographystyle{amsalpha}
%\setcounter{page}{3}
%\listoffigures

\begin{abstract}
We prove existence and uniqueness of global solutions for a class of reaction-advection-anisotropic-diffusion systems whose reaction terms have a "triangular structure". We thus extend previous results to the case of time-space dependent anisotropic diffusions and with time-space dependent advection terms. The corresponding models are in particular relevant for transport processes inside porous media and in situations in which additional migration occurs.
The proofs are based on optimal $L^p$-maximal regularity results for the general time-dependent linear operator dual to the one involved in the considered systems. As an application, we prove global well-posedness for a prototypical class of chemically reacting systems with mass-action kinetics, involving networks of reactions of the type $C_1+\ldots+C_{P-1} \rightleftharpoons C_{P} $. Finally, we analyze how a classical a priori $L^2$-estimate of the solutions, which holds with this kind of nonlinear reactive terms, extends to our general anisotropic-advection framework. It does extend with the same assumptions for isotropic diffusions and is replaced by an $L^{(N+1)/N}$-estimate in the general situation.
\end{abstract}
\maketitle
%%%%%%%%%%%%%%%%%%%%%%%%%%%%%%%%%%%%
%\tableofcontents
% \begin{equation}\label{1}
% \left\{
% \begin{array}{rcll}
%  \partial_t c_1+\mathrm{div}[-d_1(t,x)\nabla c_1 +c_1 u(t,x)]&=&-c_1c_2+c_3		 	\\[1mm]
% \partial_t c_2+\mathrm{div}[-d_2(t,x)\nabla c_1 +c_2 u(t,x)]&=&-c_1c_2+c_3		 	&\on (0,+\infty)\times \Omega\;,\\[1mm]
% \partial_t c_3+\mathrm{div}[-d_3(t,x)\nabla c_1 +c_3 u(t,x)]&=&+c_1c_2-c_3		 	\\
% \end{array}
% \right.
% \end{equation}
\section{Introduction}
%%%%%%%%%%%%%%%%%%%%%%%%%%%%%%%%%%%%%%%%%%%%%%%%%%%%%%%%%%%%%%%%%%%%%%%%%%%%%%%%%%%%%%%%%%%%%%%%%%%%%%%%%%%%%%%%%%%%%%%%%%%%%%%%%%%
%%%%%%%%%%%%%%%%%%%%%%%%%%%%%%%%%%%%%%%%%%%%%%%%%%%%%%%%%%%%%%%%%%%%%%%%%%%%%%%%%%%%%%%%%%%%%%%%%%%%%%%%%%%%%%%%%%%%%%%%%%%%%%%%%%%
Mathematical models for reacting systems are usually based on Continuum Physics, where the central role is played by the
partial mass balances
\begin{equation}
\label{partial-mass}
\partial_t \rho_i + \mathrm{div} (\rho_i v_i) = f_i,
\end{equation}
where $\rho_i$ is the mass density and $v_i$ denotes the continuum mechanical velocity of constituent $i$.
The right-hand side $f_i$ is the mass production for species $i$ due to (chemical) reactions.
The total mass balance of the reacting mixture follows by summation of the balances \eqref{partial-mass}
over all constituents, where one defines the total mass density $\rho$ and the barycentric velocity $v$ as
\begin{equation}
\label{rho-v-def}
\rho=\sum_i \rho_i, \qquad \rho v = \sum_i \rho_i v_i.
\end{equation}
Conservation of total mass corresponds to $\sum_i f_i =0$, hence the mixture obeys the continuity equation
\begin{equation}
\label{continuity}
\partial_t \rho + \mathrm{div} (\rho v) = 0.
\end{equation}
Since the barycentric velocity also appears in the momentum balance, it is common to split the
mass flux of species $i$ into an advective (also called convective) and a relative (also called molecular) contribution
according to
\begin{equation}
\label{flux-decomp}
\rho_i v_i = \rho_i v + J_i.
\end{equation}

Since chemical reactions do not change the total number of atoms, the reaction rates are more conveniently expressed for
number densities, respectively for molar concentrations $c_i:=\rho_i / M_i$, where $M_i>0$ denotes the molar mass of species $i$.
Then \eqref{partial-mass} assumes the standard form of a species balance equation, namely
\begin{equation}
\label{species-balance}
\partial_t c_i + \mathrm{div} (c_i v + j_i) = R_i,
\end{equation}
where $j_i:=J_i/M_i$ and $R_i$ is the molar production rate for component $i$.

Considering the mixture velocity $v$ to be given, the quantities $j_i$ and $R_i$ need to be modeled via constitutive,
material and composition dependent relations. Here the specific physics of the problem needs to be taken into account.
Let us briefly indicate a few prototype cases:
\begin{enumerate}
\item[(i)]
{\it Transport in solids.} Here, typically, no convection takes place and the transport is due to molecular diffusion,
possibly accompanied by drift or migration processes. In the simplest case of purely diffusive transport of
a dilute component, the molecular flux is often modeled according to Fick's law, i.e.
\begin{equation}
\label{Ficks-law}
j_i = - d_i \nabla c_i, \quad \mbox{ respectively } j_i = - d_i \nabla x_i,
\end{equation}
where $d_i>0$ is the (Fickean) diffusivity of species $i$ in the mixture and $x_i$ denotes the molar fraction,
defined as $x_i:=c_i / c_{\rm tot}$ with the total (molar) concentration $c_{\rm tot}=\sum_k c_k$.
Note that the two constitutive equations in \eqref{Ficks-law} are equivalent for
dilute components (i.e., species $i$ with $x_i \ll 1$) if the total concentration is homogeneous and constant.
In general, neither of them is applicable due to the occurrence of cross-effects and non-idealities; see
\cite{Bothe-MS} and the references therein.\\[-2ex]

\item[(ii)]
{\it Transport in porous media.} Here, in addition to transport via molecular diffusion and
possibly migration, convection often occurs due to pressure driven flow of the total mixture through the pores.
Back-effects of the transported chemical components onto the flow field are often small, in which case the velocity
field is determined by an independent set of equations, e.g.\ via Darcy's law.
The network of pores with its usually complex and partly random structure introduces so-called dispersive mixing which
appears like a diffusive transport of the same type as modeled via \eqref{Ficks-law}, but with usually much larger
dispersion coefficients as compared to molecular diffusivities.
Furthermore, due to possible anisotropy of the pore structure, the dispersive contribution to the flux will no
longer be isotropic. Therefore, the total molar mass flux is modeled as
$c_i v + j_i$ with the relative flux
\begin{equation}
\label{PM-flux}
j_i = - D_i \cdot \nabla c_i
\end{equation}
with a symmetric and positive definite matrix $D_i$. Note that $D_i$ models both dispersive and diffusive fluxes,
therefore it will typically be of the form $D_i=D_0 + d_i I$, where $D_0$ is a symmetric and positive definite matrix
modeling the dispersion in the mixture, while the diffusive contribution $d_i$ is species dependent.\\[-2ex]

\item[(iii)]
{\it Transport in fluids (liquids or gases).}
In this case, convective transport occurs because of free motion of the fluid mixture due to, e.g., natural convection, pressure driven flow or mechanical agitation. The latter is especially important in Chemical Engineering processes,
since strong mixing of initially segregated reactants is required to enable efficient chemical conversion which eventually
takes place on the molecular length scale. In case of turbulent flow fields, mixing is often modeled by dispersive transport
as mentioned above. Hence, the total flux will be of the same type as in \eqref{PM-flux}, but usually with isotropic dispersion/diffusion tensor.
Let us note in passing that possible back-effects onto the momentum balance would require a much more complicated
modeling; see \cite{BotheDreyer} and the literature cited there.\\[-2ex]

\item[(iv)]
{\it Transport by migration/drift.}
Besides molecular diffusion and macroscopic dispersion, directed transport phenomena occur in systems where external forces
are present which act differently on individual species. These forces will be balanced by the velocity-dependent friction between the constituents' molecules and the mixture, thus leading to an advective transport contribution, but with a
species dependent velocity which involves the external forces.
This is why we consider fluxes of type $c_i u_i + j_i$ with individual velocity fields $u_i$
and diffusion fluxes $j_i$ in the present paper.
An important example is the so-called ''electro-migration'',
which is caused by the (intrinsic plus externally imposed) electrical field, acting on charged particles.
See, e.g., \cite{BFPR12b},  \cite{GG96} and the references given there.
Directed motions are also relevant in ecological models, where they are modeled in the same way via advection terms,
often involving gradients of a quantity which characterizes the spatially heterogenous environment. For more about such
ecological models see, e.g., \cite{CantrellCosner}, \cite{CantrellAdvection}.
\end{enumerate}

\noindent
Despite the strong relevance of reaction-diffusion/dispersion-advection systems for applications (reacting flows, contaminant transport, electro-chemistry, population balances, etc.), only few papers address the question of global existence of solutions for general time-dependent anisotropic diffusions together with advective transport (see the brief survey below). A general theory of global existence in this framework is still missing. The present paper provides a rather general contribution in this setting in the case where nonlinear reactive terms exhibit a so-called "triangular structure". This structure is natural (and more or less necessary in this general setting) when one aims at proving global existence of {\em classical regular} solutions even for standard constant scalar diffusions and without any advection.

To be more specific, let us consider the following $2\times 2$ model example of such systems
\begin{equation}\label{examp}
\left.
\begin{array}{rcll}
 \partial_t c_1+\mathrm{div}[-D_1(t,x)\nabla c_1 +c_1 u_1(t,x)]&=&c_2-c_1h(c_1,c_2)		 	&\on (0,+\infty)\times \Omega\;,\\[1mm]
\partial_t c_2+\mathrm{div}[-D_2(t,x)\nabla c_2 +c_2 u_2(t,x)]&=&c_1+c_1h(c_1,c_2)		 	&\on (0,+\infty)\times \Omega\;,
\end{array}
\right\}
\end{equation}
where $D_1, D_2$ are positive definite matrices which are regular in $(t,x)$, $u_i, i=1,2$ are regular $R^N$-valued vectors (where $N$ denotes the space dimension) and $h:\R_+\times\R_+\to\R_+$ is a regular nonnegative function with $h(\cdot,0)\equiv0$.

Local well-posedness on some maximal interval $[0,T^*)$ for this system is well-known and global well-posedness (i.e. $T^*=+\infty$) would follow from the existence of uniform {\it a priori } bounds in $L^\infty(\Omega)$ on $[0,T^*)$ (see, e.g., \cite{LSU, amann89, amann93}). Also notice that the reaction terms have the so-called quasi-positivity structure (see {\bf (H4)} below) which guarantees that the solutions remain nonnegative. Using this nonnegativity and the fact that the nonlinear terms add up to a linear function of $c_1,c_2$, it follows from integrating the sum of the two equations of $(\ref{examp})$ that its solutions are uniformly bounded in $L^1(\Omega)$ on all finite subinterval of $[0,T^*)$. In the presence of {\em different diffusions}, even in the constant diagonal case $D_i=d_i\,I$ with $d_i\in (0,\infty)$ and without advection, it has been shown that uniform bounds in $L^1(\Omega)$ are not sufficient to prevent blow-up in $L^\infty(\Omega)$ in finite time for reaction-diffusion systems with even such 'good' properties (see \cite{PSch}). Therefore, existence of global solutions is a serious largely open question in general.

However, in the particular case of system $(\ref{examp})$, it is possible to exploit its 'triangular structure' namely: 1) the first nonlinear term is bounded above by a linear function of $c=(c_1,c_2)$; 2) the sum of the two nonlinear terms is bounded above by a linear function of $c$.

In the case of constant diagonal diffusion, without advection terms and with such a triangular structure, global existence of classical solutions has been proved in \cite{hollis,Mo, pruss} (see also Theorem 3.5 \cite{pierre10} for an up-to-date proof). A main tool is the use of $L^p$-maximal regularity results for the dual problem associated with the linear part of (\ref{examp}). Uniqueness also follows thanks to the regularity of the solutions.

In the present work, we extend the latter result to the case of anisotropic time-space dependent diffusions and with time-space dependent advection terms. A main tool is again the use of $L^p$-maximal regularity results for the dual of the linear operator involved in the system. We mainly rely on the optimal maximal regularity results proved in \cite{dhp}. We believe that they actually provide optimal global existence of bounded (and therefore regular) solutions in our setting. Note that our results have already proved to be useful in \cite{rolland12}, \cite{fischer13}, \cite{BFPR12b} for even more involved systems where advection terms are coupled with extra partial differential equations.

The existing work on global existence for reaction-diffusion-advection equations includes some contributions which address the full
system of (partial) mass, momentum and energy balance equations. In \cite{FeireislTrivisa}, the existence of a global weak solution is obtained under several assumptions.
In particular, the diffusive fluxes are modeled via Fick's law with diffusivities being
equal for all constituents. Note that the latter is required for consistency of the partial mass balances with the
continuity equation if Fick's law is employed for all species. But equal diffusivities are rarely met in
physical systems, showing that the simple approach via Fick's law is not sufficient.
In \cite{Texier}, the full set of balance equations is considered in one space dimension, again with Fickean diffusion
and equal diffusion coefficients.\\
\indent
In \cite{Kraeutle}, existence of global solutions for reaction-diffusion/dispersion-advection systems is obtained
in cases of anisotropic diffusion, but with exactly equal diffusion-advection operators for the different species. This simplifies quite a lot the analysis and, as already said, is rarely met in applications. Note also that global existence of solutions has also been considered for anisotropic diffusions in \cite{GoudonVasseur}.

Somewhat less related to the present study is \cite{HollisMorgan}, where a certain kind of
nonlinear advection terms is included in a reaction-diffusion-advection system for two species.
The influence of an additional advection term on global existence for a scalar equation
has also been studied by several authors. Since this does not touch the problems appearing for systems,
let us only refer to \cite{Levine} as an entry point into the relevant literature.

We apply our general global existence result to extensions of classical models from mass action kinetics' chemistry, where $P$ chemical species $C_1,\ldots,C_P$ are transformed according to the reaction scheme
$$\alpha_j^1C_1+\ldots+\alpha_j^PC_P\ {\large{\rightleftharpoons} }\ C_{i_j},\qquad
j\in\{1,...,R\},\; i_j\in \{1,\ldots, P\}.$$
Results are collected in Corollary \ref{example}.

Besides the question of global existence for general reaction-advection-anisotropic-diffusion systems, we also analyze in this paper the persistency of $L^2(Q_T)$-a priori space-independent estimates for the solutions. It is well-known that such $L^2(Q_T)$-estimates hold (and have proven to be very useful) in the case of scalar constant diffusions and without advection terms when the sum of the nonlinear terms is nonpositive or bounded above by a linear function (which is implied by the "triangular structure" assumed here; see {\bf (H5)}). We prove that these $L^2$-estimates are still valid in the diagonal case $D_i=d_i\,I$ and with advection terms $u_i$, {\em under the same regularity assumptions as for global existence}, namely
$$\nabla d_i, u_i\in L^\infty\left(0,T;L^r(\Omega)\right), r>\max\{2,N\},\;\; d_i\in C(\overline{Q_T}).$$
These estimates depend on the $L^2(\Omega)$-norm of the initial data.

It is unlikely that they still hold in the case of nondiagonal matrices $D_i$. On the other hand, we prove that an $L^{(N+1)/N}(Q_T)$-estimate holds for the solutions of the general system under essentially the same regularity assumptions (see Proposition \ref{LN+1} for a precise statement). They are obtained by duality of the famous $L^\infty$ (or even $C^\alpha$) Krylov-Safonov estimates for non-divergence parabolic equations with (only) bounded coefficients. It is interesting to notice that this estimate depends only on the $L^1(\Omega)^P$-norm of the initial data.
%
%
%\nd The article is organized as follows. In Section \ref{S2}, we prove global wellposedness for a regularized version of system $(\ref{eq_th_global})$. In Section \ref{S3}, we use this result to prove Theorem \ref{th_triang}. Section~\ref{Sec_ex}  is devoted to the proof of Corollary 1. Finally, Sections \ref{S5} and \ref{S6} contain respectively the proofs of Theorems \ref{thm_frl} and \ref{th:slow}.
\section{Main Results}
Let us describe in more detail the class of systems we are interested in. Throughout this work, $\Omega$ is an open bounded subset of $\R^N$, whose boundary $\Sigma:=\partial \Omega$ is of class $C^2$. For $T>0$, we write $Q_T=\Omega \times (0,T)$
and $\Sigma_T=\partial \Omega\times (0,T)$. We denote by $\nu$ the outer unit normal vector field on $\partial \Omega$ and $\partial_\nu c$ is the outer normal derivative of a function $c=(c_1, \ldots , c_P)$.
The modulus of continuity of a function $h:\overline{Q}_T\rightarrow \R$ is defined as
\begin{equation}
\label{def:modulus}
\omega_{h,T}:\R_+\rightarrow \R_+ \cup \{+\infty\},\quad
 \delta \mapsto \sup \big\{ |h(t,x)-h(\bar{t},\bar{x})|\;;\; |t-\bar{t}|+\|x-\bar{x}\|\leq \delta  \big\}.
\end{equation}
Slightly abusing this notation, we write $\omega_{h,J}$ if $t,\bar t$ in \eqref{def:modulus} are restricted to a subinterval $J\subset (0,T)$.
\vspace{2mm}
We consider the RDA-system
\begin{equation}
\label{eq_th_global}
\left.
\begin{array}{rcll}
\partial_t c_i+\mathrm{div}\big( -D_i(t,x)\nabla c_i +c_i u_i(t,x) \big) &=&f_i(t,x,c)	&\on (0,+\infty)\times \Omega,\\[1mm]
\big(-D_i(t,x)\nabla c_i +c_iu_i(t,x)\big) \cdot\nu & = & 0 			&\on (0,+\infty)\times \partial \Omega,\\[1mm]
c_i(0,\cdot)&=&c_i^0										&\on \Omega,
\end{array}
\right\}
\end{equation}
where $i\in \{1,\ldots,P\}$, $c=(c_1,\ldots ,c_P)$ and the diffusion tensor $D_i$ with $D_i (t,x) \in \R^{N\times N}$ models anisotropic diffusion.
Our assumptions on the data are the following, where $\R^P_+$ denotes $[0,+\infty)^P$:
\begin{enumerate}[\bf (H1)]
 \item  $c^0=(c_1^0,\ldots,c_P^0)\in L^\infty(\Omega,\R_+^P)$.\label{assumption_c0}\vspace{2mm}
 \item  $D_i=[d^{\,i}_{kl}]_{1\leq k,l \leq N}$ is symmetric and positive definite with  $d^{\,i}_{kl}\in C(\R_+\times\overline \Omega;\R)$ as well as $\nabla d^{\,i}_{kl}\in L^\infty_{loc}(\R_+;L^r(\Omega)^N)$ for some $r>\max \{ 2,N \}$. \label{assumption:d_i}\\[1ex]
\noindent
 Note that the symmetry assumption on the $D_i$ is in fact no restriction, since
 $\mathrm{div} (-D_i \nabla c_i) = \mathrm{div} (-D^{\rm sym}_i \nabla c_i + u_i^D c_i)$
 with $D^{\rm sym}_i=(D_i + D_i^{\sf T})/2$, $u_i^D = \mathrm{div} (D^{\rm sym}_i  -D_i)$ and $u_i^D$ has the same regularity as
 the velocity $u_i$.\vspace{2mm}
 \item   $u_i\in L^\infty_{loc}(\R_+;L^r(\Omega)^N)$, $r>\max \{ 2,N \}$. \label{assumption_u}\vspace{2mm}
 \item  $f \in C^1(\R_+\times \Omega\times \R^P,\R^P)$ is  quasi-positive, i.e.\
 $$ f_i(t,x,y)\geq 0\mbox{ whenever }(t,x,y)\in (0,+\infty)\times \Omega\times \R_+^P
 \mbox{ is such that } y_i=0. $$\label{assumption_pos}\vspace{-3mm}
 \item There exists a lower triangular invertible matrix $Q=[q_{ij}]_{1\leq i,j\leq P}$ with strictly positive diagonal entries and $b\in \R^P_+$ such that
\begin{equation*}
 \forall (t,x,y)\in \R_+\times \Omega\times \R_+^P,\quad Qf(t,x,y)\leq \Big(1+\sum_{j=1}^P y_j \Big)b.
\end{equation*}\label{assumption_tri}\\[-1ex]
\noindent
Above, we assume that $f(t,x,\cdot)$ is defined on all of $\R^P$ in order to avoid simple, but technical extension arguments. But all results remain valid for quasi-positive $f\in C^1(\R_+\times \Omega\times \R^P_+,\R^P)$.\vspace{2mm}
\item \label{assumption_last} $f$ has polynomial growth with respect to the last variable, i.e.\
$$ \forall T>0,\;\exists C, \gamma>0: \forall i, \; \forall (t,x,y)\in Q_T \times \R_+^P,\quad	 |f_i(t,x,y)|\leq C(1+|y|^\gamma).$$
\end{enumerate}
Let us note that,  by {\bf (H\ref{assumption:d_i})}, there exist
$$0<\underline d (T)\leq  \overline{d}(T)<+\infty$$
such that
\begin{equation}\label{not:d}
\underline{d}(T) |\xi |^2 \leq \langle D_i(t,x) \xi ,\xi \rangle \leq
\overline{d}(T) |\xi |^2 \; \mbox{ for all $(t,x)\in \overline{Q}_T$ and } \xi \in \R^N.
\end{equation}
%%%%%%%%%%%%%%%%%%%%%%%%%%%%%%%%%%%%%%%%%%%%%%%%%%%%%%%%%%%%%%%%%%%%%%%%%%%%%%%%%%%%%%%%%%%%%%%%%%%%%%%%%%%%%%%%%%%%%%%%%%%%%%%%%%%
\\
We prove the following well-posedness result.
\begin{theorem}\label{th_triang}
Under assumptions {\bf (H\ref{assumption_c0})-(H\ref{assumption_last})},  System $(\ref{eq_th_global})$ has a unique global nonnegative weak solution $c=(c_1,\ldots,c_P)$ in the following sense:
\begin{equation}\label{defweaksolution}
\left.
\begin{minipage}{135mm}
$\forall T>0$, $\forall i\in \{1,\ldots,P\}$, $c_i\in  C([0,T];L^2(\Omega))\cap L^{\infty}(Q_T)\cap L^2(0,T;W^{1,2}(\Omega))$;\\
$\forall \psi\in C^\infty(\overline{Q_T})$ such that $\psi(T)=0$,
\begin{equation*}
 -\int_{\Omega} c_i^0 \psi(0)+\int_{Q_T}\big(-c_i\partial_t\psi +(D_i \nabla c_i-c_i u_i)\cdot\nabla \psi\big)
 =\int_{Q_T}f_i (\cdot,\cdot, c) \psi.
\end{equation*}
\end{minipage}\right\}
\end{equation}
\nd Moreover, for any $T>0$, there exists $C>0$ depending only on $\|c^0\|_{L^\infty(\Omega)^P}$ and on
\begin{equation}\label{th2:dependence}
T,\underline d(T),\overline d(T),\omega_{d^{\,i}_{kl},T},
\|\nabla d^{\,i}_{kl}\|_{L^\infty(0,T;L^r(\Omega)^N)}, \|u_i\|_{L^\infty(0,T;L^r(\Omega)^N)}, Q, \gamma, b
\end{equation}
{\rm (}with $r>\max\{2,N\}$ from {\bf (H2), (H3)}{\rm )} such that
\begin{equation}\label{th2:est}
\|c\|_{L^\infty(Q_T)^P}+\|c\|_{L^2(0,T;W^{1,2}(\Omega)^P)}+\|\partial_t c\|_{L^2(0,T;W^{-1,2}(\Omega)^P)}\leq C.
\end{equation}
\end{theorem}
%%%%%%%%%%%%%%%%%%%%%%%%%%%%%%%%%%%%%%%%%%%%%%%%%%%%%%%%%%%%%%%%%%%%%%%%%%%%%%%%%%%%%%%%%%%%%%%%%%%%%%%%%%%%%%%%%%%%%%%%%%%%%%%%%%%
Let us briefly comment on the regularity of the data. Our purpose is to derive a global existence result under weak assumptions on the advection fields $u_i$, so that it may be used, in particular, for fixed-point arguments like in \cite{BFPR12b}.
With assumptions {\bf (H\ref{assumption:d_i})-(H\ref{assumption_u})}, we cannot expect maximal regular solutions to~$(\ref{eq_th_global})$, say solutions $c\in W^{1,2}(0,T;L^2(\Omega))\cap L^2(0,T;W^{2,2}(\Omega))$; indeed, this would require more regularity for the boundary data and for the traces
of the coefficients $D_i (t,x)$.
Our assumptions, however, are sufficient to deduce maximal regular solutions for the dual problem of~$(\ref{eq_th_global})$, at least for frozen boundary condition coefficients,
which will allow us to derive the estimate $(\ref{th2:est})$.

%%%%%%%%%%%%%%%%%%%%%%%%%%%%%%%%%%%%%%%%%%%%%%%%%%%%%%%%%%%%%%%%%%%%%%%%%%%%%%%%%%%%%%%%%%%%%%%%%%%%%%%%%%%%%%%%%%%%%%%%%%%%%%%%%%%
To prove Theorem~1, we first derive a global existence result (Proposition \ref{prop_globalexist}) under extra regularity assumptions on the data, based on local existence and continuation theorems from \cite{amann89}. This result is interesting in itself, since it provides extra regularity for the solution with smooth data. Global existence is shown by proving that any solution is {\it a priori }bounded in $L^\infty(Q_T)$ for any $T>0$. For this purpose, we first derive bounds in $L^p(Q_T)$ for any finite $p$ by a duality method, where maximal regularity theory plays a crucial role. Here, we mainly rely on \cite{dhp}.
Then, by a classical result on parabolic equations, the solutions are bounded in $L^\infty(Q_T)$. Since the $L^\infty(Q_T)$-bounds only require assumptions {\bf (H1)-(H6)}, we get the existence of weak solutions for non-smooth coefficients
by approximation of the data. Finally, we prove uniqueness for these solutions.
%%%%%%%%%%%%%%%
%%%%%%%%%%%%%%%
%%%%%%%%%%%%%%%
%%%%%%%%%%%%%%%

Next we apply this result to extensions of classical models from mass action kinetics' chemistry
to the reaction-diffusion-advection case, where we consider the following situation: $P$ chemical species $C_1,\ldots,C_P$ with molar concentrations $c_1,\ldots,c_P$ are placed in a bounded domain $\Omega$, representing an isolated reactor
(a vessel, say). For mass transport, we use the same assumptions as for Theorem~1.
The species $C_i$ are involved in $R\geq1$ chemical reactions which occur simultaneously.
For $j\in\{1,\ldots,R\}$, the $j^{th}$ chemical reaction reads
$$\alpha_j^1C_1+\ldots+\alpha_j^PC_P\ \; \underset{k_j^b}{\overset{k_j^f}{{\rightleftharpoons} }  }  \ \;\beta_j^1C_1+\ldots+\beta_j^PC_P,$$
where $ \alpha_j=(\alpha_j^1,\ldots,\alpha_j^P)$, $ \beta_j=(\beta_j^1,\ldots,\beta_j^P)\in \N_0^P$ are the so-called
stoichiometric coefficients.
For $j\in\{1,\ldots,R\}$ we define stoichiometric vectors as $\omega_j=\beta_j-\alpha_j\in \Z^P$.
The quantities $k_j^f,\, k_j^b$ are the reaction rate coefficients of the forward and backward reaction path,
respectively, and we assume them to be constant, having in mind isothermal systems.
The ratio $\kappa_j:=k_j^b / k_j^f$ is the so-called equilibrium constant for
the $j^{\rm th}$ reaction and, to see its role, we write below $k_j$ and $k_j \kappa_j$ instead of $k_j^f$ and $k_j^b$, respectively.
We also use the notation $c^\gamma:=\Pi_{i=1}^P c_i^{\gamma_i}$ for $ \gamma \in \N_0^P$ and $c=(c_1,\ldots,c_P)\in\R_+^P$. Moreover, we assume
\begin{enumerate}
\item [$\mathbf{(a_1)} $]
\label{app:ass2}
$\omega_1,\ldots,\omega_R\in\R^P$ are linearly independent (network without loops);
\item [$\mathbf{(a_2)} $]
\label{app:ass1}
$\forall j\in \{1,\ldots,R\}$, $r_j(c)=c^{\alpha_j}-\kappa_j c^{\beta_j}$ (mass action kinetics);
\item [$\mathbf{(a_3)} $]
\label{app:ass3}
$\exists e\in (0,+\infty)^{P}$ such that, for all $j\in \{1,\ldots,R\} $, $\langle e ,\omega_j\rangle=0$ (conservation of atoms);
\item [$\mathbf{(a_4)} $]
\label{app:ass1bis}
$\forall j\in \{1,\ldots,R\}$, $\beta_j$ is a permutation of $(1,0,\ldots,0)\in \N_0^P$ (single product per reaction).
\end{enumerate}
Notice that $(\mathbf{a}_4)$ requires the reactions to be of the type
$$\alpha_j^1C_1+\ldots+\alpha_j^PC_P\ {\large{\rightleftharpoons} }\ C_{i_j}\;,\; i_j\in \{1,\ldots, P\}.$$

\noindent
The equations describing the evolution of $(c_1,\ldots,c_P)$ are
\begin{equation}\label{app:main}
\left.
\begin{array}{rcll}
\partial_t c_i+\mathrm{div}\big(-D_i(t,x)\nabla c_i + u_i(t,x) c_i\big)&=&\sum_{j=1}^R \omega_j^ik_jr_j(c)	 &\on (0,+\infty)\times \Omega,\\[0.5ex]
\big( -D_i(t,x)\nabla c_i +u_i(t,x) c_i\big) \cdot\nu&=&0 					&\on (0,+\infty)\times \partial \Omega,\\[0.5ex]
c_i(0,\cdot)&=&c_i^0										&\on \Omega,
\end{array}
\right\}
\end{equation}
where $i$ runs from 1 to $P$.
\begin{corollary}\label{example}
Under assumptions {\bf (H\ref{assumption_c0})-(H\ref{assumption_u})} and
$(\mathbf{a}_{\small 1})-(\mathbf{a}_{\small 4})$,
System $(\ref{app:main})$ has a unique global solution in the sense of $(\ref{defweaksolution})$.
\end{corollary}
The proof consists in applying Theorem 1, after checking that the reaction terms satisfy assumptions {\bf (H4)-(H6)}. A technical difficulty is to check that $(\ref{app:main})$ has the ``triangular structure'' {\bf (H5)}. This requires a careful re-sorting of chemical reactions and components, and crucially relies on  $(\mathbf{a}_4)$ to get linear upper bounds.

%%%%%%%%%%%%%%%
%transition to the second part
%%%%%%%%%%%%%%%

As a special case of Corollary 1, we obtain for instance existence and uniqueness of a {\em  global regular solution} $c=(c_1,c_2,c_3):(0,+\infty)\times \Omega\rightarrow \R_+^3$ for the following system associated with the typical chemical reaction $C_1+C_2\rightleftharpoons C_3$, where $u_i$  are the advecting velocity fields. The system reads as
\begin{equation}
\label{eq_Rn}
%\tag{$R^n$}
\left.\begin{array}{rcll}
%\label{eq_c1}
\pa_tc_1+\divv(-D_1\nabla c_1+u_1\hspace{1pt} c_1)&=&-k(c_1c_2-\kappa c_3)&\on (0,\infty)\times\Omega,\vspace{1mm}\\
%\label{eq_c2}
\pa_tc_2+\divv(-D_2\nabla c_2+u_2\hspace{1pt} c_2)&=&-k(c_1c_2-\kappa c_3) &\on (0,\infty)\times\Omega,\vspace{1mm}\\
%\label{eq_c3}
\pa_tc_3+\divv(-D_3\nabla c_3+u_3\hspace{1pt} c_3)&=&\hspace{3mm}k(c_1c_2-\kappa c_3)&\on (0,\infty)\times\Omega,\vspace{1mm}\\
%\label{eq_bdy}
\big( -D_i\nabla c_i +c_iu_i\big) \cdot\nu&=&0&\on (0,\infty)\times\pa\Omega, \; i=1,2,3,\vspace{1mm}\\
%\label{eq_ic}
c_i(0,\cdot)&=&c_i^0&\on \Omega,  \; i=1,2,3.
\end{array}\right\}
\end{equation}

Besides global existence, several interesting questions have been considered in the literature for this model system. We may for instance wonder what happens to its solution when the rate constant $k$ tends to $\infty$. Considering that, here, diffusions are space-time dependent tensors and that they are perturbed by space-time dependent advection terms, this is a quite new problem. The case of constant coefficients ($D_i=d_iI, d_i\in (0,\infty)$) and without advection has been studied in \cite{BPR12}. Convergence in $L^2(Q_T)^3$ to a limit system is proved. One of the main tool is a -now  classical- $L^2(Q_T)$-estimate valid in any dimension for the solutions of these systems as soon as the nonlinearity satisfies (only) one inequality like $\sum_iq_i f_i(u)\leq 0$  with $q_i\in (0,\infty)$ for all $i$ (which is here implied by {\bf (H3)} ).

We are not going to study here the passage to the limit as $k\to\infty$ in System (\ref{eq_Rn}) (we refer to \cite{rolland12}-\cite{fischer13}) for results in this direction). But it is interesting to analyze whether or not the above mentioned $L^2(Q_T)$-estimate still holds or not. It turns out that, if the diffusions are scalar (namely $D_i=d_i\,I$ where $d_i=d_i(t,x)$ is a real-valued function), then {\em this $L^2(Q_T)$-estimate does hold exactly under the same regularity assumptions as in Theorem 1}.  We will prove the following result.

\begin{proposition}\label{propL2} Assume {\bf (H1)-(H6)} as in Theorem \ref{th_triang}. Assume moreover that, for $i=1,...,P$,  $D_i(t,x)=d_i(t,x)\,I$ where $d_i:[0,\infty)\times \overline{\Omega} \mapsto (0,\infty)$. Then, the solution of System (\ref{eq_th_global}) satisfies
$$\sum_i \|c_i\|_{L^2(Q_T)}\leq C\,[1+\sum_i \|c_i^0\|_{L^2(\Omega)}]\;for\; all\; T>0$$
with some $C$ depending only on the same quantities as in (\ref{th2:dependence}) of Theorem \ref{th_triang}.
\end{proposition}

It does not seem that this $L^2$-estimate remains valid for general matrices $D_i$. However, it is interesting to notice that an $L^{(N+1)/N}(Q_T)$-estimate does hold for the solutions of Theorem 1 in any space dimension $N$. This is obtained, by duality, from the famous Krylov-Safonov estimates for parabolic non-divergence operators. More precisely, we have

\begin{proposition}\label{LN+1} Assume {\bf (H1)-(H6)} of Theorem \ref{th_triang} with, moreover, $r>\max\{2,N+1\}$.
Then the solution of System (\ref{eq_th_global}) satisfies
$$\sum_i\|c_i\|_{L^{(N+1)/N}(Q_T)}\leq C\,[1+\sum_i \|c_i^0\|_{L^1(\Omega)}],$$
for some $C$ depending on the same quantities as in (\ref{th2:dependence}) of Theorem \ref{th_triang}.
\end{proposition}

\begin{remark}\label{remark0} In the diagonal case $D_i=d_i I$, the existence result of Theorem \ref{th_triang} remains valid for functions $d_i:[0,\infty)\times\overline{\Omega}\to (0,\infty)$ which are less regular in time and more regular in space. Instead of {\bf (H2)-(H3)}, we may assume
$$d_i\in C\big(\R_+\times\overline{\Omega};[\underline{d},\overline{d}]\big),\; \nabla d_i,u_i\in L^s_{loc}(\R_+;L^r(\Omega)^N),$$
$$2\leq s\leq+\infty,\; 2\leq N< r\leq +\infty,\; \frac{1}{s}+\frac{N}{2r}<\frac{1}{2}.$$
The proof is essentially the same as in the case $s=+\infty$ (see \cite{fischer13} for the required changes), the main point being that maximal $L^p$-regularity holds for the dual problem (\ref{eq_lemdual1}) under the above conditions as proved in \cite{dhp}. In this diagonal case, the boundary condition in the dual problem (\ref{eq_lemdual1}) is simply equivalent to $\langle \nabla  \hat\Psi, \nu \rangle =\widehat{\vartheta} \on \Sigma_T$. For general matrices $D_i$ and similar assumptions on the $d^i_{kl}$, it is not clear what to choose as a good auxiliary boundary condition in the dual problem (\ref{eq_lemdual1}) with only $L^s$-regularity in time of the $d^i_{kl}$. Thus we do not know whether the same extension can be made in general.
\end{remark}
%%%%%%%%%%%%%%%
%%%%%%%%%%%%%%%
%%%%%%%%%%%%%%%
%%%%%%%%%%%%%%%
%%%%%%%%%%%%%%%%%%%%%%%%%%%%%%%%%%%%%%%%%%%%%%%%%%%%%%%%%%%%%%%%%%%%%%%%%%%%%%%%%%%%%%%%%%%%%%%%%%%%%%%%%%%%%%%%%%%%%%%%%%%%%%%%%%%
%%%%%%%%%%%%%%%%%%%%%%%%%%%%%%%%%%%%%%%%%%%%%%%%%%%%%%%%%%%%%%%%%%%%%%%%%%%%%%%%%%%%%%%%%%%%%%%%%%%%%%%%%%%%%%%%%%%%%%%%%%%%%%%%%%%
\section{Global existence for an approximate system}\label{S2}
%%%%%%%%%%%%%%%%%%%%%%%%%%%%%%%%%%%%%%%%%%%%%%%%%%%%%%%%%%%%%%%%%%%%%%%%%%%%%%%%%%%%%%%%%%%%%%%%%%%%%%%%%%%%%%%%%%%%%%%%%%%%%%%%%%%
%%%%%%%%%%%%%%%%%%%%%%%%%%%%%%%%%%%%%%%%%%%%%%%%%%%%%%%%%%%%%%%%%%%%%%%%%%%%%%%%%%%%%%%%%%%%%%%%%%%%%%%%%%%%%%%%%%%%%%%%%%%%%%%%%%%
%
%\BR WITH THE NEW IDEAS, WE MIGHT CALL IT A THEOREM? \ER
%
\begin{proposition}\label{prop_globalexist}
In addition to {\bf (H\ref{assumption_c0})-(H\ref{assumption_last})}, assume
$$d^{\,i}_{kl}\in C^{2}(\R_+\times \overline{\Omega},\R),\;   u_i\in C^{2}(\R_+\times \overline{\Omega},\R^N),
\;c^0\in C^2(\overline \Omega,\R_+^P).$$
 Then system $(\ref{eq_th_global})$ has a unique global classical nonnegative solution
\begin{equation}\label{regular-solution}
c=(c_1,\ldots,c_P)\in C (\R_+;C(\overline \Omega)^P)\cap C^{1}((0,+\infty);C(\overline \Omega)^P)\cap C((0,+\infty);C^2(\overline \Omega)^P),
\end{equation}
and it satisfies the estimates $(\ref{th2:est})$ with a constant $C$ as characterized in \eqref{th2:dependence}.
\end{proposition}
%%%%%%%%%%%%%%%%%%%%%%%%%%%%%%%%%%%%%%%%%%%%%%%%%%%%%%%%%%%%%%%%%%%%%%%%%%%%%%%%%%%%%%%%%%%%%%%%%%%%%%%%%%%%%%%%%%%%%%%%%%%%%%%%%%%
%%%%%%%%%%%%%%%%%%%%%%%%%%%%%%%%%%%%%%%%%%%%%%%%%%%%%%%%%%%%%%%%%%%%%%%%%%%%%%%%%%%%%%%%%%%%%%%%%%%%%%%%%%%%%%%%%%%%%%%%%%%%%%%%%%%
%\item Extra regularity up to $t=0$ can be obtained provided $c^0$ satisfies appropriate compatibility conditions on $\partial\Omega$. For instance, if $c^0\in W^{s,p}(\Omega,\R_+^P)$ with $ \frac{N}{p}<s<\min(1+\frac{1}{p},2-\frac{N}{p})$, the map $(t,c^0)\mapsto c(t,\cdot),\; \R_+\times W^{s,p}(\Omega,\R_+^P) \rightarrow W^{s,p}(\Omega,\R_+^P)$ is continuous (see \cite{amann93}).
%%%%%%%%%%%%%%%%%%%%%%%%%%%%%%%%%%%%%%%%%%%%%%%%%%%%%%%%%%%%%%%%%%%%%%%%%%%%%%%%%%%%%%%%%%%%%%%%%%%%%%%%%%%%%%%%%%%%%%%%%%%%%%%%%%%

\nd As in \cite{pierre10}, the global existence is based on certain $L^p$-estimates which are obtained by duality. But in the present context we need to carefully control the trace of the solutions on the boundary of $\Omega$. At this point we therefore recall a few facts from the theory of $L^p$-maximal regularity which can be found in \cite{dhp}.
For $0\leq\tau < \tau+\delta$, we let $Q_\tau^\delta =(\tau,\tau + \delta)\times \Omega$
and $\Sigma_\tau^\delta =(\tau,\tau + \delta)\times \Sigma$.
The dual problems which shall be employed are parabolic PDEs of the type
\begin{equation}
\label{dual}
\left.\begin{array}{rcll}
\partial_t w +\div (-D(t,x)\nabla w)+u(t,x)\cdot \nabla w =f(t,x)
&\on Q_\tau^\delta,\vspace{1mm}\\
\langle D^\Sigma(t,x)\nabla w, \nu \rangle=g(t,x) &\on \Sigma_\tau^\delta,\vspace{1mm}\\
w(\tau)=w_0 &\on \Omega.
\end{array}\right\}
\end{equation}
We will only involve time-independent coefficients in the Neumann condition
in order to avoid additional regularity requirements, but recall a few facts for
the general case.
First of all, in order to apply \cite[Theorem 2.1]{dhp} to \eqref{dual}, we
bring the latter into the non-divergence form
\begin{equation}
\label{dual2}
\left.\begin{array}{rcll}
\partial_t w -D(t,x) : \nabla^2 w + u^D(t,x)\cdot \nabla w =f(t,x)
& \mbox{ in } Q_T,\vspace{1mm}\\
\langle D^\Sigma (t,x)\nabla w, \nu \rangle=g(t,x) & \mbox{ on } \Sigma_T,\vspace{1mm}\\
w(0)=w_0 & \mbox{ in } \Omega.
\end{array}\right\}
\end{equation}
Here $D:F=\sum_{k,l} d_{kl} f_{kl}$ is the double contraction between second rank tensors,
$\nabla^2 w$ denotes $[\partial_k \partial_l w]$ and $u^D:=u - \div D =u-v, v=(v_l)_l, v_l=\sum_k\partial_kd_{kl}$.
% has the same regularity as the velocity $u$.\\[2ex]

For given $1<p<\infty, \,p\neq 3$,
the PDE \eqref{dual2}, and hence \eqref{dual} as well,
has a unique strong solution (in the $L^p$-sense) if $D$ satisfies {\bf (H2)} and $D^\Sigma$ satisfies \eqref{not:d} on $\Sigma_\tau^\delta$ and $u^D$ satisfies {\bf (H3)} and the coefficients in the Neumann boundary condition satisfy
\begin{equation}
\label{reg-DN}
d^\Sigma_{kl} \in W^\kappa_s (J; L^q (\Sigma)) \cap L^s (J; W^{2\kappa}_q (\Sigma))
\end{equation}
for $s\geq 2$ and $q\geq N$ such that $\frac 2 s + \frac{N-1}{q} < 1-\frac 1 p$ and, finally, the data of the problem satisfies
$$
f\in L^p (Q_\tau^\delta),\quad
g\in W^\kappa_p (J;L^p (\Sigma))\cap L^p (J;W^{2\kappa}_p (\Sigma)),\quad
w_0 \in W^{2-\frac 2 p}_p (\Omega);
$$
here $J=(\tau,\tau + \delta)$ and $\kappa=\frac 1 2 - \frac{1}{2p}$.
In addition,  in case $p>3$, the compatibility condition $g(\tau,\cdot)=\langle D^\Sigma(\tau,\cdot)\nabla w_0, \nu \rangle$ is required.
We abbreviate the space for the Neumann data by
$$
\FF_p (\Sigma_\tau^\delta) := W^\kappa_p (J;L^p (\Sigma))\cap L^p (J;W^{2\kappa}_p (\Sigma)).
$$
Note that this is just the trace space of $\nabla w$, since $w$ lies in the maximal regularity space
$$
\WW_p (Q_\tau^\delta) := W^{1,p}(J ;L^p (\Omega))\cap L^p(J;W^2_p (\Omega)).
$$
By the results in \cite{dhp}, there exists a constant $C$ which only depends on
$\Omega$, the norms of the coefficients in the spaces given above,
on the modulus of continuity of $D$ and on the time interval $J$, such that
\begin{equation}
\label{maxreg0}
||w||_{\WW_p (Q_\tau^\delta)} \leq C \big( ||w_0||_{W^{2-2/ p}_p (\Omega)}+ ||f||_{L^p (Q_\tau^\delta)} + ||g||_{\FF_p (\Sigma_\tau^\delta)} \big).
\end{equation}
In the application to the dual problems we only need to consider the case $w_0=0$.
Precisely in this case, it turns out that a common constant $C=C(T)$ can be found for all subintervals $[\tau,\tau+\delta]$ of $[0,T]$ as it is shown later.

We are now in a position to prove the key duality estimate, for which we let
$$
\| \,|z|\, \|_{\FF_{p'} (\Sigma_\tau^\delta)'}=\sup\{\int_{\Sigma_\tau^\delta}|z|\,\sigma, \sigma\in C^\infty(\overline{Q_\tau^\delta}),\;\sigma\geq 0,\; \|\sigma\|_{\FF_{p'} (\Sigma_\tau^\delta)} \leq 1\},
$$
where $p'=p/(p-1)$. For future reference, note that for $J=[\tau,\tau + \delta]\subset [0,T]$,
the embedding estimate $\|\sigma\|_{L^{p'} (J;L^1(\Sigma))}\leq C_{\mathcal F}(T)
\|\sigma\|_{\FF_{p'} (\Sigma_\tau^\delta)}$ holds, so that
\begin{equation}\label{FTp}
\|\cdot\|_{\FF_{p'} (\Sigma_\tau^\delta)'}\leq C_{\mathcal F}(T) \|\cdot\|_{L^p (J;L^\infty(\Sigma))}.
\end{equation}
The constant $C_{\mathcal F}(T)$ depends on $T$, $\Omega$ and $p$, but not on $\tau$, $\delta$.
Observe also that, since $2\kappa < 1$,
\begin{equation}\label{abs0}
\|\,|w|\,\|_{\FF_{p'} (\Sigma_\tau^\delta)}\leq   \|w \|_{\FF_p (\Sigma_\tau^\delta)},
\end{equation}
which follows directly from the
definition of the norm in fractional Sobolev spaces of order $s\in (0,1)$.
Let us note in passing that the inequality \eqref{abs0} is, in general, strict.

Now we have the following result.
%%%%%%%%%%%%%%%%%%%%%%%%%%%%%%%%%%%%%%%%%%%%%%%%%%%%%%%%%%%%%%%%%%%%%%%%%%%%%%%%%%%%%%%%%%%%%%%%%%%%%%%%%%%%%%%%%%%%%%%%%%%%%%%%%%%
\begin{lemma}\label{lem_dualityestimate}
Given $T>0$, let $D_j=[d^j_{kl}]$ be diffusion tensors satisfying {\bf (H\ref{assumption:d_i})} and
$u_j \in C([0,T]\times \overline \Omega;\R^N)$ for $j=0,\ldots ,m$. Let
$$
v_j \in C ([0,T);C(\overline \Omega))\cap C^{1}((0,T);C(\overline \Omega))\cap C((0,T);C^2(\overline \Omega))
\; \mbox{ for } j=0,\ldots ,m
$$
satisfy
\begin{equation}\label{eq_lem_dual}
 \left.
\begin{array}{ll}
\partial_t v_0 +\div J_0 \leq \sum_{j=1}^m \Big(\alpha_j \partial_t v_j + \beta_j \div J_j + \gamma_j v_j \Big)
  &\on Q_T,\vspace{1mm}\\
\langle J_j , \nu \rangle =0 \;\mbox{ for } j=0,\ldots ,m
 &\on \Sigma_T
\end{array}\right\}
\end{equation}
with $J_j = -D_j(t,x)\nabla v_j + u_j (t,x) v_j$ and $\alpha_j, \beta_j, \gamma_j \in\R$.
Moreover, let $v_0$ satisfy $v_0 \geq 0$.
Then, for every $p>r/(r-N)$ with $r>\max\{2,N\}$ from {\bf (H2), (H3)},
there exist $\delta_0 >0$ and $C>0$, depending only on
\begin{equation}\label{lem:dual:constants}
\mbox{ $p$, $T$, $\underline d(T)$, $\overline d(T)$,
$\|\nabla d^{\,j}_{kl}\|_{L^\infty(0,T;L^r(\Omega))}$, $\|u_j\|_{L^\infty(0,T;L^r(\Omega))}$, $\omega_{d^{\,0}_{kl},T}$ and $\alpha_j, \beta_j, \gamma_j$,}
\end{equation}
such that
\begin{equation}\label{eq_lem_dual_concl}
 \|v_0\|_{L^p(Q_\tau^\delta)} + || v_0 ||_{\FF_{p'} (\Sigma_\tau^\delta)'} \leq C
\Big(  \| v_0 (\tau)\|_{L^p(\Omega)} + \sum_{j=1}^m \big( \| v_j(\tau)\|_{L^p(\Omega)} + \|v_j\|_{L^p(Q_\tau^\delta)}
+||\, | v_j| \, ||_{\FF_{p'} (\Sigma_\tau^\delta)'} \big) \Big)
\end{equation}
for every choice of $\, 0\leq \tau < \tau + \delta \leq T$ with $\delta \leq \delta_0$.
\end{lemma}
%%%%%%%%%%%%%%%%%%%%%%%%%%%%%%%%%%%%%%%%%%%%%%%%%%%%%%%%%%%%%%%%%%%%%%%%%%%%%%%%%%%%%%%%%%%%%%%%%%%%%%%%%%%%%%%%%%%%%%%%%%%%%%%%%%%

%%%%%%%%%%%%%%%%%%%%%%%%%%%%%%%%%%%%%%%%%%%%%%%%%%%%%%%%%%%%%%%%%%%%%%%%%%%%%%%%%%%%%%%%%%%%%%%%%%%%%%%%%%%%%%%%%%%%%%%%%%%%%%%%%%%
{\nd \bf Proof.}
Given $J=[\tau,\tau+\delta] \subset [0,T]$, let $\Theta\in C_0^\infty(Q_\tau^\delta)^+$,
$\vartheta\in C_0^\infty((\tau, \tau+\delta ); C^2 (\Sigma ))^+$
and consider the dual problem
\begin{equation*}
 -[\partial_t \Psi +\div (D_0 \nabla \Psi)+u_0 \cdot \nabla \Psi]=\Theta \on Q_\tau^\delta,\quad
 \langle D_0^\Sigma  \nabla \Psi, \nu \rangle =\vartheta \on \Sigma_\tau^\delta,
 \quad \Psi(\tau+\delta ,\cdot)=0 \on \Omega,
\end{equation*}
where $D_0^\Sigma (x):= D_0 (\tau + \delta ,x)$. Note that we consider the dual problem with inhomogeneous
boundary data.
Equivalently, we have $\Psi(t,\cdot)=\hat \Psi (2\tau+\delta -t,\cdot)$, where $\hat \Psi$ satisfies
\begin{equation}\label{eq_lemdual1}
 \partial_t \hat\Psi+\div (-D_0\nabla \hat\Psi)+u_0\cdot \nabla \hat\Psi=\widehat{\Theta} \on Q_\tau^\delta,\quad
 \langle D_0^\Sigma \nabla  \hat\Psi, \nu \rangle =\widehat{\vartheta} \on \Sigma_\tau^\delta,\quad
 \hat\Psi(\tau,\cdot)=0 \on \Omega,
\end{equation}
with $\Theta(t,\cdot):=\widehat{\Theta}(2\tau+\delta-t,\cdot), \vartheta(t,\cdot)= \widehat{\vartheta}(2\tau+\delta-t,\cdot)$. We are going to employ $L^{p'}$-maximal regularity for \eqref{eq_lemdual1}, where $p'=p/(p-1)$ is the dual exponent to $p$.
%Since we only want large $p$ for the primal problem, we may assume that $p'$ is as close to 1 as needed.
As mentioned above, we are going to apply \cite[Theorem 2.1]{dhp} to \eqref{eq_lemdual1},
brought into the non-divergence form
\begin{equation}\label{eq_lemdual2}
 \partial_t \hat\Psi - D_0(t,x) : \nabla^2 \hat\Psi +u(t,x) \cdot \nabla \hat\Psi=\widehat{\Theta} \on Q_\tau^\delta,\quad
 \langle D_0^\Sigma (x) \nabla  \hat\Psi, \nu \rangle =\widehat{\vartheta} \on \Sigma_\tau^\delta,\quad
 \hat\Psi(\tau,\cdot)=0 \on \Omega,
\end{equation}
where $u:=u_0 - \div D_0$ has the same properties as $u_0$.
Let us briefly comment on the assumptions required in \cite[Theorem 2.1]{dhp}.
With $\Omega$ a bounded domain with $C^2$-boundary, the result can be applied without
assumptions ''at infinity''. The assumption {\bf (H2)} implies normal ellipticity of the interior symbol as well as
the Lopatinskii-Shapiro condition at the boundary. The required regularity of the top-order coefficients and
of the advection term follow immediately from {\bf (H2), (H3)}.
In case $p' >3$ a compatibility condition is to be imposed, but which is trivially fulfilled since
$\hat \Psi (\tau ,\cdot)=0=\widehat{\vartheta}(\tau,\cdot)$.
The only remaining condition is on the boundary coefficients of the dual problem (see (\ref{reg-DN})). Since the boundary
condition is independent of time, it suffices to have $d^\Sigma_{0,kl} \in W^{1-1/p'}_q (\Sigma)$ for some
$q\geq p'$ such that $(N-1)/q <1 - 1/p' =1/p\,$; recall the meaning of $N$ here, i.e.\ $\Omega \subset \R^N$.
By the assumptions on $D_0$ and trace theorems for Sobolev functions,
it holds that $d^\Sigma_{0,kl} \in W^{1-1/r}_r (\Sigma)$ with $r$ from {\bf (H2)}.
By Sobolev embedding, this implies $d^{\Sigma}_{0,kl}\in W^{1-1/p'}_q(\Sigma)$ for any $q$
such that $1-\frac N r > 1 - \frac 1 {p'} - \frac {N-1} {q} = \frac 1 p - \frac {N-1} {q}$.
Since only the condition $1-\frac N r >  \frac 1 p $ remains in the limit case $q\to \infty$,
all conditions are satisfied (choosing a sufficiently large $q$) in case $p' < r/N$ and the latter holds
by our assumption on $p$.
Hence \cite[Theorem 2.1]{dhp} applies, showing that the dual problem \eqref{eq_lemdual1}
has $L^{p'}$-maximal regularity, i.e.
\begin{equation}\label{D1}
\|\Psi\|_{\WW_{p'} (Q_\tau^\delta)}
\leq C(T)
\big( \|\Theta\|_{L^{p'}(Q_\tau^\delta)}+\|\vartheta\|_{\FF_{p'} (\Sigma_\tau^\delta)} \big)
\end{equation}
for any data $\Theta$ and $\vartheta$ from the corresponding data spaces.
The estimate \eqref{D1} is uniform with respect to $\tau$ and $\delta$, i.e.\ $C(T)$ is
independent of $\tau$ and $\delta$. The latter follows from the fact the functions $\widehat{\Theta}$
and $\widehat{\vartheta}$ have compact support in time in $(\tau, \tau + \delta)$, hence can be
extended by zero to all of $[0,T]$ (yielding extensions $\tilde \Theta$, $\tilde \vartheta$)
without changing their norms in the respective spaces.
Indeed, let $\tilde \Psi$ be the unique strong solution of
\begin{equation}\label{eq_lemdual3}
\partial_t \tilde\Psi - D_0(t,x) : \nabla^2 \tilde\Psi +u(t,x) \cdot \nabla \tilde\Psi=\tilde\Theta \on Q_T,\quad
\langle D_0^\Sigma (x) \nabla  \tilde\Psi, \nu \rangle =\tilde\vartheta \on \Sigma_T,\quad
\tilde\Psi(0,\cdot)=0 \on \Omega
\end{equation}
and note that $\tilde \Psi_{|Q_\tau^\delta} = \hat \Psi$ as well as
$\|\Psi\|_{\WW_{p'} (Q_\tau^\delta)} = \| \hat\Psi\|_{\WW_{p'} (Q_\tau^\delta)}$.
Then \eqref{D1} follows from
\begin{equation}
\| \hat\Psi\|_{\WW_{p'} (Q_\tau^\delta)}\leq\| \tilde\Psi\|_{\WW_{p'} (Q_T)}
\leq C(T)
\big( \|\tilde \Theta\|_{L^{p'}(Q_T)}+\|\tilde \vartheta\|_{\FF_{p'} (\Sigma_T)} \big)
= C(T)
\big( \|\Theta\|_{L^{p'}(Q_\tau^\delta)}+\|\vartheta\|_{\FF_{p'} (\Sigma_\tau^\delta)} \big).
\end{equation}

Having estimate \eqref{D1} at hand, we now test the primary problem (\ref{eq_lem_dual}) with $\Psi$, where
it is important to note that $\Psi \geq 0$ due to $\Theta \geq 0$ and $\vartheta \geq 0$.
Multiplying (\ref{eq_lem_dual}) by $\Psi$ and integrating over $Q_\tau^\delta$,
partial integration leads to the following inequality (where $\sum_j$ means $\sum_{j=1}^{m}$)
\begin{equation}\label{grandestim}
\left.
\begin{array}{l}
\int_{Q_\tau^\delta}v_0 \Theta+\int_{\Sigma_\tau^\delta} v_0 \vartheta
\leq\int_\Omega \big( v_0(\tau)-\sum_j\alpha_j v_j(\tau) \big) \Psi(\tau)\\[1.5ex]
- \int_{Q_\tau^\delta} \sum_jv_j \big( \alpha_j \partial_t\Psi
+\beta_j (\div (D_j(t,x) \nabla\Psi)+u_j(t,x) \cdot \nabla\Psi )- \gamma_j \Psi \big)\\[1.5ex]
+\sum_j\beta_j\int_{\Sigma_\tau^\delta} v_j \langle D_j(t,x) \nabla\Psi, \nu \rangle
+\int_{\Sigma_\tau^\delta} v_0 \langle (D_0^\Sigma(x)-D_0(t,x)) \nabla\Psi, \nu\rangle.
\end{array}
\right\}
\end{equation}
We now employ the estimate
$$
\left|\int_{\Sigma_\tau^\delta} v_0 \langle (D_0^\Sigma(x)-D_0(t,x))\nabla\Psi, \nu\rangle \right|
\leq \omega_{D_0,J}  \int_{\Sigma_\tau^\delta} |v_0|\,|\nabla\Psi|\leq
\omega_{D_0,J} \| v_0 \|_{\FF_{p'} (\Sigma_\tau^\delta)'}
\sum_q\|\,|\partial_{x_q} \Psi|\,\|_{\FF_{p'} (\Sigma_\tau^\delta)},
$$
where $\omega_{D_0,J}:=\max \{\omega_{d^{\,0}_{kl},J}: k,l=1,\ldots N \}$; recall
the definition of $D_0^\Sigma$ and the fact that $v_0\geq 0$.
Using \eqref{abs0}, i.e.
\begin{equation}\label{abs}
\|\,|\partial_{x_q}\Psi|\,\|_{\FF_{p'} (\Sigma_\tau^\delta)}\leq
\|\partial_{x_q}\Psi\|_{\FF_{p'} (\Sigma_\tau^\delta)},
\end{equation}
which holds since $2\kappa:=1-1/p' < 1$,
we deduce that
\begin{equation}\label{FT}
\left.
\begin{array}{l}
\int_{Q_\tau^\delta}v_0 \Theta+\int_{\Sigma_\tau^\delta} v_0 \vartheta\leq \|v_0(\tau)-\sum_j\alpha_j v_j(\tau)\|_{p}
\|\Psi(\tau)\|_{p'}
+C_0(\sum_j\|v_j\|_{L^p(Q_\tau^\delta)} )\|\Psi\|_{\WW_{p'} (Q_\tau^\delta)}\\[1ex]
+C_0\,\left[\sum_j\||v_j|\|_{\FF_{p'} (\Sigma_\tau^\delta)'}+\omega_{D_0,J}
\| v_0 \|_{\FF_{p'} (\Sigma_\tau^\delta)'}\right] \sum_q\|\partial_{x_q}\Psi\|_{\FF_{p'} (\Sigma_\tau^\delta)}
\end{array}
\right\}
\end{equation}
with a constant $C_0$ depending only on the data. Employing (\ref{D1}) and
$$\sum_q\|\partial_{x_q}\Psi\|_{\FF_{p'} (\Sigma_\tau^\delta)}
\leq C(T) \|\Psi\|_{\WW_{p'} (Q_\tau^\delta)} $$
by Sobolev imbedding (with a constant which only depends on $p'$, $T$ and $\Omega$),
inequality \eqref{FT} implies
$$\int_{Q_\tau^\delta}v_0 \Theta+\int_{\Sigma_\tau^\delta} v_0 \vartheta\leq
C(T)\,\left[K_\tau^\delta +C_0 \, \omega_{D_0,J} \|v_0 \|_{\FF_{p'} (\Sigma_\tau^\delta)'}\right] \,\left[\|\Theta\|_{L^{p'}(Q_\tau^\delta)}+\|\vartheta\|_{\FF_{p'} (\Sigma_\tau^\delta)}\right]$$
with the abbreviation
$$K_\tau^\delta:= \|v_0 (\tau)-\sum_j\alpha_j v_j (\tau)\|_{p}+C_0\sum_j\|v_j\|_{L^p(Q_\tau^\delta)}+
C_0\sum_j\||v_j|\|_{\FF_{p'} (\Sigma_\tau^\delta)'}.$$
Since $\Theta\in C_0^\infty(Q_\tau^\delta)^+$ and
$\vartheta\in C_0^\infty((\tau, \tau+\delta ); C^2 (\Sigma ))^+$
are arbitrary and $v_0\geq 0$, we deduce that
\begin{equation}\label{dualestimate}
\|v_0 \|_{L^p(Q_\tau^\delta)}+\| \, v_0 \, \|_{\FF_{p'} (\Sigma_\tau^\delta)'}\leq
C(T) \,\left[K_\tau^\delta+C_0 \, \omega_{D_0,J} \| \, v_0 \, \|_{\FF_{p'} (\Sigma_\tau^\delta)'}\right].
\end{equation}
We now choose $\delta>0$ so small that $\omega_{D_0,J} \leq (2 C(T) C_0 )^{-1}$ and obtain
\begin{equation}\label{mainlp}
\|v_0 \|_{L^p(Q_\tau^\delta)}+\| \, v_0 \, \|_{\FF_{p'} (\Sigma_\tau^\delta)'}\leq 2 C(T)\,K_\tau^\delta,
\end{equation}
hence \eqref{eq_lem_dual_concl}.\vspace{0.05in}
\cqfd \vspace{0.05in}
%%%%%%%%%%%%%%%%%%%%%%%%%%%%%%%%%%%%%%%%%%%%%%%%%%%%%%%%%%%%%%%%%%%%%%%%%%%%%%%%%%%%%%%%%%%%%%%%%%%%%%%%%%%%%%%%%%%%%%%%%%%%%%%%%%%

%%%%%%%%%%%%%%%%%%%%%%%%%%%%%%%%%%%%%%%%%%%%%%%%%%%%%%%%%%%%%%%%%%%%%%%%%%%%%%%%%%%%%%%%%%%%%%%%%%%%%%%%%%%%%%%%%%%%%%%%%%%%%%%%%%%
{\nd \bf Proof of Proposition \ref{prop_globalexist}.}
Existence of a unique nonnegative regular
(i.e., with the regularity as stated in \eqref{regular-solution})
solution $c$ of \eqref{eq_th_global} on a maximal time interval $[0,T^*)$, $0<T^*\leq +\infty$ is a consequence of Theorem 15.1 in \cite{amann93} (see also \cite{amann89}).

To prove $T^*=+\infty$, it is sufficient to prove that, if $T^*<+\infty$, then $c\in L^\infty(Q_{T^*})^P$ (see \cite[Theorem 3]{amann89}). We will actually prove that, given $T\leq T^*, T<+\infty$, there exists $\delta$ depending only on $T$ and the data such that for all $[\tau,\tau+\delta)\subset [0,T)$, $\|c\|_{L^\infty([\tau, \tau+\delta))^P}\leq C$ with $C$ depending only on $T$ and the data as well. As explained below, it will follow at the same time that $T^*=+\infty$ and that the solution $c$ satisfies the estimate $\|c\|_{L^\infty(Q_T)^P}\leq C$ in  (\ref{th2:est}).

Let $c=(c_j)_{1\leq j\leq P}$ be the solution on $[0,T^*)$. We fix $T\leq T^*, T<\infty$. Let us set $W:=\sum_{j=1}^P c_j$. We will estimate $W$ on intervals $[\tau,\tau+\delta) \subset [0,T)$. Recall the notation $Q_\tau^\delta=(\tau,\tau+\delta)\times\Omega, \Sigma_\tau^\delta=(\tau,\tau+\delta)\times \partial\Omega$. Using $Q= [q_{ij} ]$ and $b=(b_1,\ldots,b_P) $ defined in {\bf (H\ref{assumption_tri})}, we introduce the solution $z_i$ of
\begin{equation}
\label{eq:z}
 \left.
\begin{array}{rcl}
 q_{ii}\partial_t z_i +\mathrm{div} (-D_i \nabla z_i + u_i z_i)=(1+W)b_i&\on& Q_\tau^\delta,\\[1mm]
(-D_i\nabla z_i +z_i u_i) \cdot\nu=0 &\on& \Sigma_\tau^\delta,\\[1mm]
z_i(\tau,\cdot)=0&\on& \Omega.
\end{array}
\right\}
\end{equation}
Let $p>(N+2)/2$. We then have
\begin{equation}\label{zinfty}
\|z_i\|_{L^\infty(Q_{\tau}^{(t-\tau)})}\leq C[1+\|W\|_{L^p(Q_{\tau}^{(t-\tau)})}]
\;\mbox{ for all } t\in [\tau, T).
\end{equation}
This follows from Theorem III.7.1  in \cite{LSU} for Dirichlet boundary conditions. But the same result holds for Neumann boundary conditions as well. The proof is carried out in detail  for isotropic diffusions in \cite{fischer13, rolland12} and in the appendix of \cite{BotheRolland14} and the proof given there remains valid without any changes
for anisotropic diffusion tensors of the type considered here.

Notice that the right-hand side in \eqref{eq:z} is linear in $(1+W)$, whence the corresponding
linear dependence in the estimate just given. It is important to note that the constant $C$ in \eqref{zinfty}
does not depend on $\tau, t$ but only on $T$.

Below, we employ the abbreviation
$\mathcal{E}_j c_j= \partial_t c_j +\mathrm{div} (-D_j \nabla c_j + u_j c_j)$
and, analogously,  $\mathcal{E}_j z_j$ for $z_j$ instead of $c_j$.
By \eqref{eq_th_global} and {\bf (H5)}, we have
\begin{align}
 q_{ii}\mathcal{E}_i c_i & =q_{ii}f_i(t,x,c)	\leq (1+W)b_i-\sum_{j=1}^{i-1}q_{ij}f_j(t,x,c)=(1+W)b_i-\sum_{j=1}^{i-1}q_{ij}\mathcal{E}_j c_j \label{eq_pr_thgl1}
\end{align}
for every $i\in \{1,\ldots,P\}$.
Employing \eqref{eq:z}, inequality $(\ref{eq_pr_thgl1})$ can be rewritten as
\begin{equation}
\label{inequ:ci}
 q_{ii}\mathcal{E}_i c_i \leq  q_{ii}\mathcal{E}_i z_i
 -\sum_{j=1}^{i-1}q_{ij} \mathcal{E}_j c_j.
\end{equation}
We apply Lemma~\ref{lem_dualityestimate} to \eqref{inequ:ci} with $m=i$, $v_0=q_{ii} c_i$ and, say,
$v_j=-c_j$, $\alpha_j=\beta_j=q_{ij}$ for $j=1, \ldots , i-1$ and $v_i=z_i$, $\alpha_i=\beta_i=q_{ii}$.
Denoting
$$
\|\cdot\|_{\PP_{\tau,t}}
=\|\cdot\|_{L^p((\tau,t)\times\Omega)}+\|\cdot\|_{\FF_{p'}(\Sigma_{\tau}^{(t-\tau)})'}
\;\mbox{ for }t\in [\tau, \tau+\delta]
$$
and choosing $p>1$ above so large that also $p>r/(r-N)$ with $r>N$ from {\bf (H2), (H3)} holds, Lemma~\ref{lem_dualityestimate} yields
$$\|c_i\|_{\PP_{\tau,t}}\leq C \Big( \sum_{j=1}^i \|c_j(\tau)\|_p +\sum_{j<i}\|c_j\|_{\PP_{\tau,t}}+\|z_i\|_{\PP_{\tau,t}}\Big)$$
with a constant $C$ depending on $T$ but not on $\tau , t$.
%
%
%The value of $\delta$ depends only on $T,p'$ and on the regularity of the $D_1, D_2, u_1, u_2$ (and the $\theta_i, \Omega$). We choose one $\delta$ which fits to all operators $\zeta\to \EE_j\zeta:=\partial_t\zeta+div(-D_j\nabla\zeta+u_j\zeta)$.\\
%
%Going back to our situation with the $c_i$, with the notation of our paper, we have $$q_{ii}\EE_ic_i\leq -\sum_{j<i}q_{ij}\EE_jc_j+\EE_iz_i,$$ where the $z_i$ are defined as before, except that we work only on $[\tau,T)$ and we choose $z_i(\tau)=0$.
%
By induction, we obtain
\begin{equation}\label{cilp}
\|c_i\|_{\PP_{\tau,t}}\leq C
\Big( \sum_{j=1}^i \|c_j(\tau)\|_p +\sum_{j\leq i}\|z_j\|_{\PP_{\tau,t}} \Big).
\end{equation}
Summing over $i$ from 1 to $P$, we deduce
\begin{equation}\label{wallp}
\|W\|_{\PP_{\tau,t}}\leq C
\sum_{j=1}^P \big( \|c_j(\tau)\|_p + \|z_j\|_{\PP_{\tau,t}} \big).
\end{equation}
Combining this inequality with \eqref{zinfty} and using (see (\ref{FTp}) )
\begin{equation}\label{FT'}
\|\cdot\|_{\FF_{p'} (\Sigma_\tau^\delta)'}\leq C_{\FF}(T) \|\cdot\|_{L^p\left(\tau,\tau+\delta;L^\infty(\Sigma)\right)},
\end{equation}
we get
\begin{equation}\label{zz}
\sum_{j=1}^P \|z_j\|_{L^\infty(Q_{\tau}^{(t-\tau)})}\leq C \Big[ 1+
\sum_{j=1}^P \big(\|c_j(\tau) \|_p + \|z_j\|_{L^p(Q_{\tau}^{(t-\tau)})}
+\|z_j\|_{L^p (\tau,\tau+\delta;L^\infty(\Sigma))} \big) \Big].
\end{equation}
Taking the $p$-th power and employing the fact that
$\|\phi\|^p_{L^p(\Omega)}+\|\phi\|^p_{L^\infty(\Sigma)}\leq C \|\phi\|^p_{L^\infty(\Omega)}$ for continuous functions $\phi$ on $\overline{\Omega}$, we deduce
$$\forall t\in (\tau,\tau+\delta),\;\; \sum_{j=1}^P\|z_j(t)\|^p_{L^\infty(\Omega)}\leq
C[1+\sum_{j=1}^P \|c_j(\tau) \|_p +\sum_{j=1}^P\int_\tau^t\|z_j(s)\|^p_{L^\infty(\Omega)}ds].$$
This is a Gronwall inequality for $[\tau, \tau+\delta ]\ni t \mapsto \sum_{j=1}^P \|z_j(t)\|^p_{L^\infty(\Omega)}$.
Going back to (\ref{wallp}), we obtain
\begin{equation}\label{cilpfinal}
\max_i\|c_i\|_{L_p(Q_\tau^\delta)}\leq C \big( \max_i\|c_i(\tau)\|_p+1 \big).
\end{equation}
with a constant $C$ which depends on $T$ but not on $\tau, \delta$.

We now restrict again the choice of $p$ so that $p >\gamma  (N+2)/2$ with $\gamma$ defined in the growth condition {\bf (H6)} on $f$. It follows that, with various constants $C$ depending on the same quantities,
$$\max_i \|f_i(c)\|_{L^p(Q_\tau^\delta)}\leq C\left[1+\|c\|^\gamma_{L^{p\gamma}(Q_\tau^\delta)^P}\right]\leq C[1+\max_i\|c_i(\tau)\|^\gamma_{L^{p\gamma}(\Omega)}].$$
We again exploit Theorem III.7.1  in \cite{LSU} (in the appropriately modified version for Neumann conditions as explained for the solution $z_i$ in the comments following (\ref{eq:z}) and (\ref{zinfty})) to obtain
\begin{equation}
\label{ineq:tau-delta}
\max_i\|c_i\|_{L^\infty(Q_\tau^\delta)}\leq  \Phi \big( \max_i\|c_i(\tau)\|_\infty +1 \big)
\end{equation}
with some (nonlinear) function $\Phi$. Note that the same $\Phi$ applies for any $[\tau,\tau +\delta)\subset [0,T)$ with
$\delta \in (0,\delta_0 ]$, where $\delta_0 >0$ comes
from application of Lemma~\ref{lem_dualityestimate} above.

We next apply \eqref{ineq:tau-delta} successively to the subintervals $[k\delta ,k \delta +\delta]$ for $k=0,\ldots ,n-1$, where $\delta =T/n$ with sufficiently large $n\in \N$ so that $\delta \leq \delta_0$,
the latter coming from Lemma~\ref{lem_dualityestimate}.
This yields the {\it a priori} estimate
\begin{equation}
\label{ineq:0-T}
\max_i\|c_i\|_{L^\infty(Q_T)}\leq \Phi^n \Big( \max_i\|c_i^0 \|_\infty +1\Big).
\end{equation}

As explained at the beginning of this proof, this implies at the same time $T^*=+\infty$ and the estimate $\|c\|_{L^\infty(Q_T)}\leq C$ contained in (\ref{th2:est}), where $C$ depends on the quantities listed in Proposition \ref{prop_globalexist}.

To get the other estimates in (\ref{th2:est}), we multiply the equation (\ref{eq_th_global}) in $c_i$ by $c_i$ and integrate over $Q_T$ to obtain
$$\frac{1}{2}\|c_i(T)\|^2_{L^2(\Omega)}+\underline{d}(T)\|\nabla c_i\|^2_{L^2(Q_T)}\leq \frac{1}{2}\|c_i^0\|^2_{L^2(\Omega)}+\int_{Q_T}|c_iu_i\nabla c_i|+c_i|f_i(c)|$$
$$\leq \frac{1}{2}\|c_i^0\|^2_{L^2(\Omega)}+\frac{\underline{d}(T)}{2}\|\nabla c_i\|^2_{L^2(Q_T)}+C\int_{Q_T}|c_iu_i|^2+c_i|f_i(c)|.$$
Using $u_i\in L^\infty(0,T;L^2(\Omega)), c\in L^\infty(Q_T)^P$, we deduce $\nabla c_i\in L^2(Q_T)$ with associated bounds depending on the same quantities. Finally, we go back to equation (\ref{eq_th_global}) to obtain the estimate of $\partial_t c_i$ in the space $L^2(0,T;W^{-1,2}(\Omega))$.\cqfd

\begin{remark}\label{remarkun}{\rm
If one had $b=0$ in {\bf (H5)}, we would have $z_i=0$ in (\ref{eq:z}). Then, estimate (\ref{cilp}) would directly give estimate (\ref{cilpfinal}) and the main part of the proof of Proposition \ref{prop_globalexist} would then be quite simpler.

If the boundary trace of $D_0$ had the regularity \eqref{reg-DN} with $p$ replaced by $p'$,
then we could apply Lemma~\ref{lem_dualityestimate} directly with $\tau=0$ and $\delta=T$ and prove Proposition~\ref{prop_globalexist} also more directly.
But note that {\bf (H2)} does not impose any time regularity of the $\nabla d^0_{kl}$. Therefore,
we have to use an auxiliary boundary conditions for the dual problem and this requires a partition of $[0,T]$ into subintervals in the proof of a priori $L^\infty$-bounds in
Proposition \ref{prop_globalexist}.
}

As already noticed in Remark \ref{remark0}, if $D_i=d_i$, where $d_i:\R_+\times\overline{\Omega}\to (0,\infty)$, then the boundary condition in (\ref{eq_lemdual1}) is simply equivalent to $\langle \nabla  \hat\Psi, \nu \rangle =\widehat{\vartheta} \on \Sigma_T$.  Thus the term $\int_{\Sigma_\tau^\delta} v_0 \langle (D_0^\Sigma(x)-D_0(t,x)) \nabla\Psi, \nu\rangle$ in (\ref{grandestim}) is equal to zero and we do not need to restrict the size of $\delta$ anymore. The $L^\infty$-estimate may then be proven directly on $[\tau, \tau+\delta)=[0,T)$.
\end{remark}

%%%%%%%%%%%%%%%%%%%%%%%%%%%%%%%%%%%%%%%%%%%%%%%%%%%%%%%%%%%%%%%%%%%%%%%%%%%%%%%%%%%%%%%%%%%%%%%%%%%%%%%%%%%%%%%%%%%%%%%%%%%%%%%%%%%

%%%%%%%%%%%%%%%%%%%%%%%%%%%%%%%%%%%%%%%%%%%%%%%%%%%%%%%%%%%%%%%%%%%%%%%%%%%%%%%%%%%%%%%%%%%%%%%%%%%%%%%%%%%%%%%%%%%%%%%%%%%%%%%%%%%
% \begin{remark}
% In \cite{amann93}, the results we used from Chapters 14, 15 and 16 are stated for quasilinear time-independent operators. To see that they are still valid for the time-dependent case, it is sufficient to ``artificially'' add $t$ in equations $(\ref{eq_th_global})$ replacing $c,u$ and $f$ by $\tilde c=(c_1,\ldots,c_P,t)$, $\tilde u=(u_1,\ldots,u_P,0)$ and $\tilde f=(f_1,\ldots,f_P,1)$.
% \end{remark}
%
%%%%%%%%%%%%%%%%%%%%%%%%%%%%%%%%%%%%%%%%%%%%%%%%%%%%%%%%%%%%%%%%%%%%%%%%%%%%%%%%%%%%%%%%%%%%%%%%%%%%%%%%%%%%%%%%%%%%%%%%%%%%%%%%%%%
\section{Proof of Theorem \ref{th_triang}}\label{S3}
{\it Existence. }
%%%%%%%%%%%%%%%%%%%%%%%%%%%%%%%%%%%%%%%%%%%%%%%%%%%
Let $T>0$. We approximate $c^0_i, D_i=[d^{\,i}_{kl}]$ and $u_i$ in System $(\ref{eq_th_global})$ by smooth functions $c^{0n}_i, D_i^n =[d^{\,i,n}_{kl}],u^n_i$ such that
\begin{equation*}
 c^{0n}_i  \underset{n\rightarrow +\infty}{\longrightarrow} c_i^0 \mbox{ in }L^2(\Omega),\;
d^{\,i, n}_{kl} \underset{n\rightarrow +\infty}{\longrightarrow} d^{\,i}_{kl} \mbox{ in }L^2(Q_T),\; u^n_i\underset{n\rightarrow +\infty}{\longrightarrow}u \mbox{ in }L^2(Q_T)^N
\end{equation*}
and such that $(c^{0n}_i)_{n\in \N}$ is bounded in $L^\infty(\Omega)^+$,
$\omega_{d^{\,i, n}_{kl},T}\leq \omega_{d^{\,i}_{kl},T}$,
$\underline{d} |\xi |^2 \leq \langle D_i(t,x) \xi ,\xi \rangle \leq \overline d |\xi |^2$,
$(u^n_i)_{n\in\N}$ and $(\nabla d^{\,i,n}_{kl})_{n\in\N}$ are bounded in $L^\infty(0,T;L^r(\Omega)^N)$.
According to Proposition $\ref{prop_globalexist}$, System $(\ref{eq_th_global})$ with data $( c^{0n}_i,D_i^n,u^n_i)$ has a unique solution $c^n: [0,T]\times \Omega\rightarrow \R_+^P$. Moreover, estimate $(\ref{th2:est})$ is satisfied and guarantees that $\|c^n\|_{L^\infty(Q_T)^P}$, $\|c^n\|_{L^2(0,T;W^{1,2}(\Omega)^P)}$ and $\|\partial_t c^n\|_{L^2(0,T;W^{-1,2}(\Omega)^P)}$ are bounded, independently of $n$. By Corollary~4 in \cite{simon86}, $(c^n)_{n\in \N}$ is relatively compact in $L^2(Q_T)^P$ and therefore has a subsequence that converges a.e. in $Q_T$.\saut
Let $T_k\nearrow \infty$. We denote by ${c^n_i}_{|[0,T_k]}$ the restriction of $c^n_i$ to $Q_{T_k}$. Using the above results, there exists $c:(0,+\infty)\times \Omega\rightarrow \R_+^P$ such that, up to a diagonal extraction, we have, as $n\rightarrow +\infty$:

$\quad\forall i\in \{1,\ldots,P\},\;\forall k\in \N,$
\begin{equation}\label{compactness}
\left.
\begin{array}{rcll}
{c^n_i}_{|[0,T_k]}&\longrightarrow& c_i &\mbox{ in }L^p(Q_{T_k})\mbox{ for any }p<+\infty \mbox{ and a.e. };\\[1mm]
f_i(t,x,c^n)_{|[0,T_k]}&\longrightarrow& f_i(t,x,c) &\mbox{ in }L^p(Q_{T_k})\mbox{ for any }p<+\infty\;;\\[1mm]
\nabla c^n_{i\;|[0,T_k]}&\longrightarrow& \nabla c_i &\mbox{ weakly in }L^2(Q_{T_k})^N;\\[1mm]
\partial_t c_i^n &\longrightarrow& \partial_t c_i&\mbox{ weakly in }L^2(0,T_{k};W^{-1,2}(\Omega)).
\end{array}
\right\}
\end{equation}
As $c^n$ is a classical solution of $(\ref{eq_th_global})$, for all $T>0$ and all $\psi\in C^\infty(\overline{Q_T})$ with $\psi(T)=0$, we have
\begin{equation}\label{eq_weak_approx}
-\int_\Omega c_i^{0n}\psi(0)+\int_{Q_T}\big(-c_i^n\pa_t\psi+(D_i^n\nabla c_i^n-c_i^nu^n_i)\cdot\nabla\psi\big) =\int_{Q_T}f_i(t,x,c^n)\psi.
\end{equation}
Using $(\ref{compactness})$, we can pass to the limit as $n\rightarrow +\infty$ in $(\ref{eq_weak_approx})$ for any $T>0$,
so that $c$ satisfies \eqref{defweaksolution} up to $c\in C([0,T];L^2 (\Omega))$; the latter can be shown as follows.
Since $c_i\in L^\infty(Q_T)\cap L^2((0,T);W^{1,2}(\Omega))$ and $\partial_tc_i\in L^2(0,T;W^{-1,2}(\Omega))$, for any $\zeta \in C_0^\infty(\Omega)$, $c_i\zeta \in L^2((0,T);W^{1,2}_0(\Omega))\cap W^{1,2}(0,T;W^{-1,2}(\Omega))\hookrightarrow C([0,T];L^2(\Omega))$ (see, e.g., \cite{Evans}) and, consequently, $c_i\in C([0,T];L^2_{loc}(\Omega))\cap L^\infty(Q_T)\hookrightarrow C([0,T];L^2(\Omega))$.\\[0.5ex]
Estimates $(\ref{th2:est})$, which we already proved for smooth solutions in Proposition \ref{prop_globalexist}, are inherited by the weak solutions since the norms can only decrease when passing to weak limits.
\vspace{3mm}\\
%%%%%%%%%%%%%%%%%%%%%%%%%%%%%%%%%%%%%%%%%%%%%%%%%%%
%%%%%%%%%%%%%%%%%%%%%%%%%%%%%%%%%%%%%%%%%%%%%%%%%%%
{\it Uniqueness. }
%%%%%%%%%%%%%%%%%%%%%%%%%%%%%%%%%%%%%%%%%%%%%%%%%%%
Let $T>0$ and let $c,\hat{c}$ be two solutions of $ (\ref{defweaksolution})$ on $Q_T$ with the same initial data, and set $w_i:=c_i-\hat{c_i}$. Below, $C>0$ denotes a generic constant depending only on $T$ and the data of $(\ref{eq_th_global})$.\saut
Let us first use a formal computation to show that $w_i=0$, and justify it afterwards. We have
\begin{equation}\label{eq_triang0}
\left.
\begin{array}{lll}
 \partial_t w_i + \textrm{div}(- D_i\nabla w_i +w_iu_i )=f_i(c)-f_i(\hat c) \on Q_T,\\[0.5ex]
\big( - D_i \nabla w_i +w_i u_i \big) \cdot \nu=0 \on \Sigma_T, \;\; w_i(0)=0 \on \Omega.
\end{array}
\right\}
\end{equation}
Letting $t\in (0,T)$, multiplying $(\ref{eq_triang0})$ by $w_i$ and integrating by parts on $Q_{t}$, we get
\begin{equation}\label{eq_triang6}
 \frac{1}{2} \int_\Omega w_i(t)^2 +\int_{Q_{t}} \langle D_i \nabla w_i, \nabla w_i \rangle
 =\int_{Q_{t}} w_i u_i\cdot \nabla w_i +\int_{Q_{t}} [f_i(c)-f_i(\hat c)] w_i.
\end{equation}
Since $c, \hat c\in L^\infty(Q_T)$ and $f_i$ is locally Lipschitz continuous,
\begin{equation}\label{eq_triang5}
 \exists C=C(T)>0: \int_{Q_{t}} [f_i(c)-f_i(\hat c)] w_i \leq C  \int_{Q_{t}} \Big(\sum_{j=1}^P |w_j|\Big )|w_i|\leq C  \sum_{i=1}^P\int_{Q_{t}}  w_i^2.
\end{equation}
Due to $u_i\in L^\infty(0,T;L^r(\Omega)^N)$, we have
(with $r^*>2$ such that $\frac{1}{r^*}+\frac{1}{r}+\frac{1}{2}=1$)
\begin{align}
 \int_{Q_{t}} w_i u_i\nabla w_i &\leq \int^{t}_0 \|w_i\|_{L^{r^*}(\Omega)} \|u_i\|_{L^r(\Omega)^N} \|\nabla w_i\|_{L^2(\Omega)^N}\nonumber\\
				&\leq C \int^{t}_0 \|w_i\|_{L^{r^*}(\Omega)} \|\nabla w_i\|_{L^2(\Omega)^N}\nonumber\\
				&\leq  \eps \|\nabla w_i\|_{L^{2}(Q_{t})^N}^2 + C_\eps \| w_i\|_{L^2(Q_{t})}^2,\label{eq_triang4}
\end{align}
where $\eps>0$ is arbitrarily small.
In \eqref{eq_triang4}, we employed the fact that
\begin{equation*}
 \forall \eps>0,\; \exists C_\eps >0:\forall w \in L^{r^*}(\Omega),\; \|w\|_{L^{r^*}(\Omega)}\leq \eps \|\nabla w\|_{L^2(\Omega)^N}+C_\eps \|w\|_{L^2(\Omega)}.
\end{equation*}
The latter follows from the compact embedding of $W^{1,2}(\Omega)$ into $L^{r^*}(\Omega)$, which holds since $r>N$ implies $-\frac{1}{2}+\frac{1}{N}>-\frac{1}{r^*}$.
Using inequalities $(\ref{eq_triang5})$ and $(\ref{eq_triang4})$ in $(\ref{eq_triang6})$,
exploiting $\langle D_i \xi ,\xi \rangle \geq \underline d |\xi |^2$,
summing over $i$ and choosing $\eps$ small enough, we get
\begin{equation}\label{eq_triang7}
  \frac{1}{2} \sum_{i=1}^P \int_\Omega w_i(t)^2 +\frac{\underline d}{2}\sum_{i=1}^P \int_{Q_{t}}|\nabla w_i|^2 \leq C\sum_{i=1}^P \int_{Q_{t}} w_i^2
  \quad \mbox{ for all } t\in(0,T).
\end{equation}
Then Gronwall's lemma yields $w_i=0$ for all $i$, i.e.\ $c=\hat c$.\vspace{1mm}
Let us briefly indicate how this computation can be justified for weak solutions.
It is clear that $(\ref{eq_triang4})$ still holds, i.e.\ we only need to justify \eqref{eq_triang6}.
The starting point is
\begin{equation}\label{psi}
\int_{Q_T}-w_i\partial_t\psi+(D_i\nabla w_i-w_iu_i)\cdot\nabla\psi=
\int_{Q_T}[f_i(c)-f_i(\hat{c})]\psi,
\end{equation}
valid for all $\psi\in C^\infty(\overline{Q_T})$ with $\psi(T)=0$. By density, and thanks to
$$w_i\in L^2(0,T:H^1(\Omega))\cap L^\infty(Q_T),\;
D_i\in L^\infty(Q_T),\; u_i\in L^\infty(0,T;L^r(\Omega)),\; r>\max\{2,N\},$$
this remains valid for all $\psi\in L^2(0,T;H^1(\Omega))$ with $\partial_t\psi\in L^1(Q_T)$ and $\psi(T)=0$. Given $t\in (0,T)$, to obtain \eqref{eq_triang6}, we would like to choose $\psi(s)=\chi_{(0,t)}(s)w_i(s)$ in (\ref{psi}). Since this $\psi$ is not regular enough in time, we rather first choose $\psi_h(s)=\frac{1}{2h}\int_{s-h}^{s+h}\psi(\sigma)d\sigma$ with $h>0$ and then let $h\to 0+$. Note that $\psi_h, \nabla \psi_h$ converge in $L^2(Q_T)$ to $\psi,\nabla\psi$ which easily allows to pass to the limit in the integrals where these terms are involved. Now, by an obvious change of time variable, we have
{\small $$\int_{Q_T}w_i\partial_t\psi_h=\int_0^t\int_\Omega w_i(s)\frac{w_i(s+h)-w_i(s-h)}{2h}=\frac{1}{2h}\int_\Omega\left[\int_h^{t+h}w_i(s-h)w_i(s)ds-\int_0^tw_i(s)w_i(s-h)ds\right],$$}
that is
$$\int_{Q_T}w_i\partial_t\psi_h=\frac{1}{2h}\int_\Omega\left[\int _t^{t+h}w_i(s-h)w_i(s)ds-\int_0^hw_i(s)w_i(s-h)ds\right].$$
Employing $w_i\in C([0,T];L^2(\Omega))$ and $w_i(0)=0$, we deduce that this term converges to $\frac{1}{2}\int_\Omega w_i(t)^2,$ which finally yields \eqref{eq_triang6}.
\cqfd
\section{Proof of Corollary  \ref{example}}\label{Sec_ex}
The proof mainly consists in reordering the reactions and the chemical components and is based on the following elementary result on matrices. We denote by $\MM_{P,R}(\R)$ the space of matrices with real entries having $P$ columns and $R$ rows, and write $\MM_{P}(\R)=\MM_{P,P}(\R)$.
\begin{lemma}\label{lem:cor}
Let $M\in \MM_{P,R}(\R)$,
$$ M= \left(
\begin{array}{c|c|c}
\;&&\\&&
\\
\hspace{-1mm}\omega_1	&\; \ldots\; 	& \omega_R\hspace{-1mm}
\\
\;&&\\&&
\end{array}
\right)\;;\quad \omega_j =
\left(
\begin{array}{c}
\hspace{-1mm}\omega_j^1 \hspace{-1mm}\\
\vspace{-1mm}
 \\ \vdots\\
\vspace{-1mm}
\\ \hspace{-1mm}\omega_j^P\hspace{-1mm}
\end{array}
\right)\in \R^P,
$$
and assume
\begin{enumerate}[$\quad (i)$]
\item $\omega_1,\ldots,\omega_R\in\R^P$ are linearly independent.\label{lem:ass2}\vspace{0.6mm}
\item $\forall j\in \{1,\ldots,R\}$, there exists a unique $i\in \{1,\ldots,P\}$ such that  $\omega_j^i>0$.\label{lem:ass1}
\item $\exists e\in (0,+\infty)^{P}$ such that $\langle e ,\omega_j\rangle=0$ for all $j\in \{1,\ldots,R\}$.\label{lem:ass3}
\end{enumerate}
Then, up to a permutation of its columns and rows,
\begin{equation}\label{lem:cor:concl}
%%%%%%%%%%%%%%%%%%%%%%%%%%%%%%%%
M={\scriptsize{
\left(\hspace{-2mm} \begin{array}{cccc}
%%%%%%%%%%%%%%%%%%%%%%%%%%%%%%%%
\begin{minipage}[t]{11mm}
\centering{
$\begin{array}[t]{|c|}
\hline
 \;\;N_1\;\;\\
\hline
\end{array}$}
\end{minipage}
%%%%%%%%%%%%%%%%%%%%%%%%%%%%%%%%
&\hspace{-2mm}
\begin{minipage}[t]{11mm}
\vspace{4mm}
\centering{
$\begin{array}[t]{|c|}
\hline
 \;\;N_2\;\;\\
\hline
\end{array}$}
\end{minipage}
%%%%%%%%%%%%%%%%%%%%%%%%%%%%%%%%
&\hspace{-2mm}
\begin{minipage}[t]{11mm}
\vspace{8mm}
\centering{ \
$\begin{array}[t]{c}
\hspace{8mm}\\
 \;\;\ddots\;\;\vspace{9mm}\\ \hspace{-21mm}
\mbox{ \huge{$\ast$}}
\end{array}$}
\end{minipage}
%%%%%%%%%%%%%%%%%%%%%%%%%%%%%%%%
&\hspace{-2mm}
\begin{minipage}[t]{11mm}
\vspace{18mm}
\centering{
\vspace{-12mm}\hspace{-8mm}\mbox{ \Large{$0$} } \hspace{8mm}\vspace{12mm}\\
$\begin{array}[t]{|c|}
\hline
 \;\;N_k\;\;\\
\hline
\end{array}$}
\vspace{15mm}   		
\end{minipage}
\end{array}\hspace{-3mm}\right)
}}\;,
\end{equation}
where $N_i$ are rows having (strictly) negative entries.
\end{lemma}
{\nd \bf Proof. }First, observe that properties $(\ref{lem:ass2})-(\ref{lem:ass3})$ are unchanged when permuting the rows or columns of $M$. We prove Lemma \ref{lem:cor} by induction on $P$. Since the vectors $e,\omega_1,\ldots,\omega_R$ are linearly independent, we have $P>R$. Using assumption $(\ref{lem:ass1})$, $M$ has exactly $R$ positive entries, so there exists $N\in \MM_{P-R,R}(\R)$ with nonpositive entries and $M_0\in \MM_R(\R)$ such that, up to a permutation of the rows,
$$
M=\left(\vspace{-1mm}
\begin{array}{ccc}
\begin{minipage}[t]{11mm}
\centering{\vspace{-3mm}
$\begin{array}[t]{|c|}
\hline
 \;\;\;N\;\;\;\\
\hline
\end{array}$}		
\end{minipage} \\
%%%%%%%%%%%%%%%%%%%%%%%%%%%%%%%%
\begin{minipage}[t]{11mm}
\centering{\vspace{-2mm}
$\begin{array}[t]{|c|}
\hline
\\
 \;\;M_0\hspace{1mm}\;\\
\;\\
\hline
\end{array}$}
\end{minipage}
\end{array}\hspace{1mm}
\right).
$$
If $N=0$, assumption $(\ref{lem:ass2})$ implies that $M_0$ is invertible, but assumption $(\ref{lem:ass3})$ yields
a nonzero vector with is orthogonal to the column vectors of $M_0$. This contradiction shows that $N\neq 0$.
Let $L_1$ be a nonzero row of $N$: by a permutation of the rows, we put $L_1$ at the top of $M$, and by a permutation of the columns, we put the strictly negative entries of $L_1$ at the top left corner, so that
$$
M=\left(\vspace{-1mm}
\begin{array}{ccc}
\begin{minipage}[t]{11mm}
\centering{\vspace{-3mm}
$\begin{array}[t]{|c|}
\hline
 \;\;\;N_1\;\;\;\\
\hline
\end{array}$}		
\end{minipage}
&\hspace{2mm}
%%%%%%%%%%%%%%%%%%%%%%%%%%%%%%%%
\begin{minipage}[t]{14mm}
\centering{\vspace{-3mm}
$\begin{array}[t]{c}
 \hspace{6mm}0\hspace{6mm}\\
\end{array}$}		
\end{minipage}
\\
%%%%%%%%%%%%%%%%%%%%%%%%%%%%%%%%
\begin{minipage}[t]{11mm}
\centering{
\vspace{4mm}
\mbox{\Large{$\;\;\;\ast$}} }
\end{minipage}
&
\begin{minipage}[t]{14mm}
\centering{\vspace{-2mm}
$\begin{array}[t]{|c|}
\hline
\\
\; \hspace{3mm}M_1\hspace{5mm}\\
\;\\
\hline
\end{array}$}
\end{minipage}
\end{array}\hspace{2mm}
\right),
$$
where $N_1$ is nontrivial and has negative entries. It is easy to check that $M_1$ satisfies assumptions $(\ref{lem:ass2})-(\ref{lem:ass3})$, hence Lemma \ref{lem:cor} holds by induction on $P$.
\cqfd\vspace{4mm}

{\nd\bf Proof of Corollary \ref{example}. }Let $M \in \MM_{P,R}(\R)$ be the matrix whose columns are $\omega_1,\ldots,\omega_R$. Using assumptions $(\mathrm{a}_{1})-(\mathrm{a}_{4})$, $M$ satisfies $({\ref{lem:ass1}})-({\ref{lem:ass3}})$ from Lemma \ref{lem:cor}. Consequently, up to a permutation of the chemical species and of the chemical reactions (which corresponds to a permutation of the rows and the columns of $M$, respectively), we can assume that $M$ has the structure as indicated in $(\ref{lem:cor:concl})$. To prove that there exists a lower triangular invertible matrix $Q\in \MM_P(\R)$ with nonnegative diagonal entries such that $QM$ has nonpositive entries, we may proceed as follows: first, recall that the multiplication of $M$ by such a matrix $Q$ corresponds to adding to each row of $M$ a positive linear combination of the above rows. We may define the matrix $Q$ as the product $Q_1\cdots Q_k$, where $Q_i$ are lower triangular invertible matrices with nonnegative entries,
satisfying
\begin{enumerate}[$\quad 1)$]
 \item The columns of $ Q_k M$ corresponding to the block $N_k$ are nonpositive. This is obtained by choosing a matrix $Q_k$ which corresponds to adding convenient positive factors of the $k^{th}$ row to the rows below.
 \item The columns of $Q_{k-1}Q_k M $ corresponding to the block $N_{k-1}$ are nonpositive. This is obtained by choosing a matrix $Q_{k-1}$ which corresponds to adding convenient positive factors of the ${(k-1)}^{th}$ row to the rows below. The crucial point is that this operation leaves the columns corresponding to $N_k$ unchanged.
 \item We iterate this procedure to build a sequence of matrices $Q_1,\ldots, Q_k$ such that the product $Q_1\cdots Q_k M$ has nonpositive entries.
\end{enumerate}
Then, denoting by $F=(F_1,\ldots,F_P)$ the reaction term in $(\ref{app:main})$, observe that
$$
\left(
\begin{array}{c}
\hspace{-1mm}F_1\hspace{-1mm}\\
\vspace{-1mm}
 \\ \vdots\\
\vspace{-1mm}
\\ \hspace{-1mm}F_P\hspace{-1mm}
\end{array}
\right)
=
\left(
\begin{array}{c|c|c}
\;&&\\&&
\\
\hspace{-1mm}\omega_1	&\; \ldots\; 	& \omega_R\hspace{-1mm}
\\
\;&&\\&&
\end{array}
\right)\cdot
\left(
\begin{array}{c}
\hspace{-1mm}k_1r_1\hspace{-1mm}\\
\vspace{-1mm}
 \\ \vdots\\
\vspace{-1mm}
\\ \hspace{-1mm}k_Pr_P\hspace{-1mm}
\end{array}
\right)
\; ; \qquad
QF=QM \left(
\begin{array}{c}
\hspace{-1mm}k_1r_1\hspace{-1mm}\\
\vspace{-1mm}
 \\ \vdots\\
\vspace{-1mm}
\\ \hspace{-1mm}k_Pr_P\hspace{-1mm}
\end{array}
\right).
$$
Using assumptions $(\mathrm{a}_{1})-(\mathrm{a}_{4})$, we have
$$ \forall j\in \{1,\ldots,P\},\quad -r_j(c)\leq \max_{j=1,\ldots,R} \kappa_j \,\sum_{i=1}^Pc_i.$$
Hence, since $QM$ has nonpositive entries, assumption {\bf (H\ref{assumption_tri})} from Theorem \ref{th_triang} is satisfied. Consequently, Theorem \ref{th_triang} applies to system $(\ref{app:main})$ and Corollary \ref{example} follows.\cqfd
%
%%%%%%%%%%%%%%%%%%%%%%%%%%%%%%%%%%%%%%%%%%%%%%%%%%%%%%%%%%%%%%%%%%%%%%%%%%%%%%%%%%%%%%%%%%%%%%%%%%%%%%%%%%%%%%%%%%%%%%%%%%%%%%%%%%%
%%%%%%%%%%%%%%%%%%%%%%%%%%%%%%%%%%%%%%%%%%%%%%%%%%%%%%%%%%%%%%%%%%%%%%%%%%%%%%%%%%%%%%%%%%%%%%%%%%%%%%%%%%%%%%%%%%%%%%%%%%%%%%%%%%%
\section{Proof of Proposition  \ref{propL2}} \label{S5}
%%%%%%%%%%%%%%%%%%%%%%%%%%%%%%%%%%%%%%%%%%%%%%%%%%%%%%%%%%%%%%%%%%%%%%%%%%%%%%%%%%%%%%%%%%%%%%%%%%%%%%%%%%%%%%%%%%%%%%%%%%%%%%%%%%%
%%%%%%%%%%%%%%%%%%%%%%%%%%%%%%%%%%%%%%%%%%%%%%%%%%%%%%%%%%%%%%%%%%%%%%%%%%%%%%%%%%%%%%%%%%%%%%%%%%%%%%%%%%%%%%%%%%%%%%%%%%%%%%%%%%%
%
%%%%%%%%%%%%%%%%%%%%%%%%%%%%%%%%%%%%%%%%%%%%%%%%%%%%%%%%%%%%%%%%%%%%%%%%%%%%%%%%%%%%%%%%%%%%%%%%%%%%%%%%%%%%%%%%%%%%%%%%%%%%%%%%%%%
\noindent
Let us first show that {\bf (H5)} implies that there exist $q_i\in (0,\infty), i=1,...,P,$ and $b_0\in \R_+$ such that
\begin{equation}\label{sumqi}
\sum_i q_i f_i(t,x,y)\leq \Big(1+\sum_{j=1}^Py_j\Big) b_0,\;\; \forall\,(t,x,y)\in \R_+\times\Omega\times \R_+^P.
\end{equation}
Indeed, we multiply the $i^{\rm th}$ line of the inequality {\bf (H5)}, namely
$$\sum_{j=1}^iq_{ij}f_j(t,x,y)\leq \Big(1+\sum_{j=1}^Py_j\Big) b_i,$$
by $\epsilon^{i-1}$ and sum over $i=1,\ldots ,P$ to obtain
$$\sum_{j=1}^P\Big[\sum_{i=j}^P\epsilon^{i-1}q_{ij}\Big]f_j(t,x,y)\leq \Big(1+\sum_{j=1}^Py_j\Big)\Big(\sum_{i=1}^{P}\epsilon^{i-1}b_i\Big).$$
Since $q_{jj}>0$ for all $j=1,\ldots ,P$,  we immediately check that, for $\epsilon>0$ small enough, we have
$$q_i:=\sum_{i=j}^{P}\epsilon^{i-1}q_{ij}=\epsilon^{j-1}\Big[q_{jj}+\epsilon\sum_{i=j+1}^{P}\epsilon^{i-j-1}q_{ij}\Big]>0
\; \mbox{ for all } i=1,\ldots ,P.$$
Whence (\ref{sumqi}) with $q_i$ defined as above and $b_0:=\sum_{i=1}^{P}\epsilon^{i-1}b_i$.

Note that we may multiply each equation (\ref{eq_th_global}) by $q_i$. We then obtain a similar system,
where $c_i, f_i$ are replaced by $\tilde{c}_i:=q_ic_i, \tilde{f}_i:=q_i f_i$ with now $\sum_i \tilde{f}_i(t,x,\widetilde{c})\leq [1+\sum_i\widetilde{c}_i/q_i]\,b_0$ by (\ref{sumqi}). Thus, without loss of generality and up to changing the value of $b_0$, we may (and will) assume that
\begin{equation}\label{tildef}
\sum_i f_i(t,x,y)\leq [1+\sum_i y_i]\,b_0,\; \forall \,(t,x,y)\in \R_+\times\Omega\times\R_+^P.
\end{equation}

At this point, recall that $D_i=d_i I$. Consider the approximate global (and regular) solution $c$ given by Proposition \ref{prop_globalexist} and let
$$W:=1+\sum_{i=1}^Pc_i,\; A:=\frac{1+\sum_{i=1}^Pd_ic_i}{1+\sum_{i=1}^Pc_i}, \;u:=\sum_{i=1}^P\frac{c_i}{W}[\nabla d_i+u_i].$$
We sum all equations (\ref{eq_th_global}) over $i=1,\ldots ,P$ and use (\ref{tildef}) to obtain
$$\partial_t \Big(\sum_i c_i\Big) +{\rm div}\Big(-\sum_id_i\nabla c_i+\sum_ic_iu_i\Big)=\sum_if_i\leq [1+\sum_i c_i]b_0$$
or
$$\partial_t \Big(\sum_i c_i\Big) +{\rm div}\Big(-\nabla [\sum_id_ic_i]+\sum_ic_i[\nabla d_i+u_i]\Big)\leq [1+\sum_i c_i]b_0,$$
or
\begin{equation}\label{eqW}
\partial_t W+{\rm div}\Big(-\nabla(A\,W)+u\,W\Big)\leq W\,b_0.
\end{equation}
Similarly, we have at the boundary
\begin{equation}\label{bdyW}
\big(-\nabla(A\,W)+u\,W\big)\cdot\nu=0 \mbox{ on } (0,+\infty)\times\partial\Omega.
\end{equation}
We have the following estimates, where $\underline{d}(T), \overline{d}(T)$ are defined in (\ref{not:d}) and
$r>\max\{2,N\}$ is from {\bf (H2), (H3)}:
$$\min\{1, \underline{d}(T)\}\leq A\leq \max\{1,\overline{d}(T)\},\;A\in C(\overline{Q_T}).$$
$$\|u\|_{L^\infty(0,T;L^r(\Omega)^N)}\leq \max_i\|\frac{c_i}{W}\|_{L^\infty(Q_T)}\sum_i\|\nabla d_i+u_i\|_{L^\infty(0,T;L^r(\Omega))}\leq \sum_i\|\nabla d_i+u_i\|_{L^\infty(0,T;L^r(\Omega))}.$$
Now the $L^2(Q_T)$-estimate of Proposition \ref{propL2} is a consequence of Lemma \ref{lem_L2_est} below applied to the solution $W$ of (\ref{eqW}-\ref{bdyW}). Indeed, Lemma \ref{lem_L2_est} says that
$$\forall j=1,\ldots ,P,\; \|c_j\|_{L^2(Q_T)}\leq \|W\|_{L^2(Q_T)}\leq C\|W(0)\|_{L^2(\Omega)}\leq C[1+\sum_i\|c_i^0\|_{L^2(\Omega)}],$$
where the constant $C$ depends on the quantities listed in Proposition \ref{propL2}. This estimate remains valid for the solution of Theorem \ref{th_triang} itself. \cqfd

% -\int_\Omega \psi(0)W_j^n(0)+\int_{Q_T} -W_j^n\pa_t\psi+\nabla (A_j^n W_j^n-W_j^nu)\cdot\nabla \psi=0,
% \end{equation}
% % for all $\psi\in C^\infty(\overline{Q_T})$ such that $\psi(T)=0$.

%%%%%%%%%%%%%%%%%%%%%%%%%%%%%%%%%%%%%%%%%%%%%%%%%%%%%%%%%%%%%%%%%%%%%%%%%%%%%%%%%%%%%%%%%%%%%%%%%%%%%%%%%%%%%%%%%%%%%%%%%%%%%%%%%%%
\begin{lemma}
\label{lem_L2_est}
Let $A\in C(\overline{Q_T})$ be such that $0<\underline a\leq A\leq \overline a <+\infty$, $u\in L^\infty(0,T;L^r(\Omega)^N)$ with $r>\max \{ 2,N \}$, $H\in L^2(Q_T)$, $b_0\in [0,\infty)$. Let $W$ be a classical solution of
\begin{equation}\label{eq_W_lem}
\Big\{
\begin{array}{ll}
\pa_tW+\divv[-\nabla(A W)+Wu]\leq H +b_0W		\on Q_T,\\[0.5ex]
-\nabla(AW)\cdot \nu+Wu\cdot\nu=0 \on \Sigma_T, \;\;W(0,\cdot)=W^0 \on \Omega.
\end{array}
\Big.
\end{equation}
Then there exists $C>0$, depending only on $T,\underline a,\overline a$ and $\|u\|_{L^\infty(0,T;L^r(\Omega)^N)}$, such that
\begin{equation}\label{eq_lem_L2_concl}
\|W^+\|_{L^2(Q_T)}\leq C\Big(\|W^0\|_{L^2(\Omega)}+\|H\|_{L^2(Q_T)}\Big).
\end{equation}
\end{lemma}
\begin{remark} {\rm
This lemma is interesting in itself. It generalizes a well-known $L^2(Q_T)$-estimate for parabolic equations (see \cite{pierre10}) to the case of variable diffusivities and additional advection term. Interestingly, the same assumption $u\in L^\infty(0,T;L^r(\Omega)^N)$ with $r>\max \{ 2,N \}$ as for Theorem 1 is required for Lemma \ref{lem_L2_est} to be valid.}

Note also that the proof of the $L^2$-estimate of Proposition \ref{propL2} requires only the structure (\ref{tildef}) (or (\ref{sumqi})) on the $f_i$ and not the full triangular structure {\bf (H5)}.
\end{remark}

%%%%%%%%%%%%%%%%%%%%%%%%%%%%%%%%%%%%%%%%%%%%%%%%%%%%%%%%%%%%%%%%%%%%%%%%%%%%%%%%%%%%%%%%%%%%%%%%%%%%%%%%%%%%%%%%%%%%%%%%%%%%%%%%%%%

%%%%%%%%%%%%%%%%%%%%%%%%%%%%%%%%%%%%%%%%%%%%%%%%%%%%%%%%%%%%%%%%%%%%%%%%%%%%%%%%%%%%%%%%%%%%%%%%%%%%%%%%%%%%%%%%%%%%%%%%%%%%%%%%%%%
{\nd \bf Proof of Lemma \ref{lem_L2_est}.}
Without loss of generality, we may assume $b_0=0$ in Lemma \ref{lem_L2_est}. Indeed, in\-tro\-ducing $\widetilde{W}:=e^{-tb_0}W$, the inequality (\ref{eq_W_lem}) may be rewritten as
$$\partial_t \widetilde{W}+\divv[-\nabla(A \widetilde{W})+\widetilde{W}u]\leq He^{-tb_0}.$$

Thus, let us assume $b_0=0$ in the following. \ER Given $\Theta\in C^\infty_0(Q_T,\R_+)$, we consider the dual problem
\begin{equation}
 \label{eq_dual_prob_2}
-[\pa_t\Psi+A\Delta\Psi+u\cdot\nabla\Psi]=\Theta \on Q_T, \;\; \pa_\nu\Psi=0 \on \Sigma_T,
\;\;\Psi(T,\cdot)=0 \on \Omega.
\end{equation}
According to Theorem~2.1 in \cite{dhp}, problem $(\ref{eq_dual_prob_2})$ has a strong solution $\Psi\geq 0$. We multiply (\ref{eq_W_lem}) by $\Psi$ and integrate over $Q_T$ and by parts to obtain
\[
\int_{Q_T} (\pa_tW+\divv(-\nabla(AW)+Wu))\Psi \leq \int_{Q_T}H\Psi ,
\]
hence
\[
-\int_{Q_T}W(\pa_t\Psi+A\Delta\Psi+u\cdot\nabla\Psi) \leq \int_{Q_T}H\Psi +\int_\Omega W^0\Psi(0).
\]
Thus
\begin{equation}
\label{eq_goal}
\Big|\int_{Q_T}W\Theta\Big|\leq\|W^0\|_{L^2(\Omega)}\|\Psi(0)\|_{L^2(\Omega)}+\|H\|_{L^2(Q_T)}\|\Psi\|_{L^2(Q_T)}.
\end{equation}
We next estimate $\|\Psi(0)\|_{L^2(\Omega)}$ and  $\|\Psi\|_{L^2(Q_T)}$ in terms of $\|\Theta\|_{L^2(Q_T)}$. Multiplying (\ref{eq_dual_prob_2}) by $-\Delta\Psi$ and integrating over $\Omega$ and by parts, using the homogeneous Neumann boundary conditions, we get
\begin{align}
-\frac12\Dt\|&\nabla\Psi\|_{L^2(\Omega)}^2+\|\sqrt{A}\Delta \Psi\|_{L^2(\Omega)}^2=-\int_\Omega u\cdot\nabla\Psi \;\Delta\Psi -\int_\Omega \Theta\Delta\Psi \nonumber\\
	&\leq \|u\|_{L^r(\Omega)^N}\|\nabla\Psi\|_{L^p(\Omega)}\|\Delta\Psi\|_{L^2(\Omega)}+\|\Theta\|_{L^2(\Omega)}\|\Delta\Psi\|_{L^2(\Omega)}\nonumber,
\end{align}
where we used H\"older's inequality and $p>1$ is defined by $1/r+1/p=1/2$.
By the Gagliardo-Nirenberg inequality (see, e.g., \cite{LSU}), there exists $C>0$ such that
\begin{equation}
\label{eq_GN}
\|\nabla\Psi\|_{L^p(\Omega)}\leq C\|\nabla\Psi\|_{W^{1,2}(\Omega)}^{N/r}\|\nabla\Psi\|_{L^2(\Omega)}^{1-N/r}.
\end{equation}
Since $u\in L^\infty(0,T;L^r(\Omega)^N)$, using Young's inequality and (\ref{eq_GN}), for every $\eps>0$
there are $C_i=C_i^\ep>0$ such that
\begin{equation*}
\begin{split}
-\frac12\Dt\|&\nabla\Psi\|_{L^2(\Omega)}^2+\underline{a}\|\Delta \Psi\|_{L^2(\Omega)}^2\\
&\leq\ep\|\Delta\Psi\|_{L^2(\Omega)}^2+C_{1}(\|\nabla\Psi\|_{L^p(\Omega)}^2+\|\Theta\|_{L^2(\Omega)}^2)\\
&\leq\ep(\|\Delta\Psi\|_{L^2(\Omega)}^2+\|\nabla \Psi\|_{W^{1,2}(\Omega)}^2) +C_{2}(\|\nabla\Psi\|_{L^2(\Omega)}^2+\|\Theta\|_{L^2(\Omega)}^2)\\
&\leq \ep\|\Delta\Psi\|_{L^2(\Omega)}^2+C_{3}(\|\nabla\Psi\|_{L^2(\Omega)}^2+\|\Theta\|_{L^2(\Omega)}^2),
\end{split}
\end{equation*}
where we used $\|\nabla\Psi\|_{W^{1,2}(\Omega)}\leq C(\Omega)\|\Delta\Psi\|_{L^2(\Omega)}$ for the last inequality. Since $\pa_\nu \Psi=0 \on \pa \Omega$ and $\Omega$ is smooth, this is a consequence of elliptic regularity; see, e.g., \cite{Brezis}. Thus, if we choose $\ep<\underline a / 2$ and apply Gronwall's lemma, using $\Psi(T)=0$, we get
\begin{equation}
\label{eq_est_Psi_1}
\sup_{0\leq t\leq T}\|\nabla\Psi(t)\|_{L^2(\Omega)}\leq C\|\Theta\|_{L^2(Q_T)} \mbox{ and } \|\Delta \Psi\|_{L^2(Q_T)}\leq C\|\Theta\|_{L^2(Q_T)}.
\end{equation}
Since $A\leq \overline{a}$, we also have $\|A\Delta \Psi\|_{L^2(Q_T)}\leq C\|\Theta\|_{L^2(Q_T)}$. Then, integration of $(\ref{eq_dual_prob_2})$ on $Q_t$ for any $t\in(0,T)$ yields
\begin{equation}\label{psilinfl1}
 \|\Psi\|_{L^\infty(0,T;L^1(\Omega))}\leq C \|\Theta\|_{L^2(Q_T)}.
\end{equation}
Finally, combining $(\ref{eq_est_Psi_1})$, $(\ref{psilinfl1})$ and using the Poincar\'e-Wirtinger inequality, we get
\[\|\Psi(0)\|_{L^2(\Omega)}+\|\Psi\|_{L^2(Q_T)}\leq C\|\Theta\|_{L^2(Q_T)},\]
whence $(\ref{eq_lem_L2_concl})$ by duality.\cqfd
%%%%%%%%%%%%%%%%%%%%%%%%%%%%%%%%%%%%%%%%%%%%%%%%%%%%%%%%%%%%%%%%%%%%%%%%%%%%%%%%%%%%%%%%%%%%%%%%%%%%%%%%%%%%%%%%%%%%%%%%%%%%%%%%%%%
%%%%%%%%%%%%%%%%%%%%%%%%%%%%%%%%%%%%%%%%%%%%%%%%%%%%%%%%%%%%%%%%%%%%%%%%%%%%%%%%%%%%%%%%%%%%%%%%%%%%%%%%%%%%%%%%%%%%%%%%%%%%%%%%%%%
%
%
%%%%%%%%%%%%%%%%%%%%%%%%%%%%%%%%%%%%%%%%%%%%%%%%%%%%%%%%%%%%%%%%%%%%%%%%%%%%%%%%%%%%%%%%%%%%%%%%%%%%%%%%%%%%%%%%%%%%%%%%%%%%%%%%%%%
\section{Proof of Proposition \ref{LN+1}}
We consider again the solution $c$ given by Proposition \ref{prop_globalexist}. It satisfies equations (\ref{eq_th_global}) namely
$$\partial_t c_i+\div \Big(-[\sum_ld^i_{kl}\partial_lc_i]_k+u_ic_i\Big)=f_i.$$
As in the proof of Proposition \ref{propL2} we will, without loss of generality, assume that (\ref{tildef}) holds.
We let $W:=1+\sum_{i=1}^P c_i$ and sum all $i$-equations to obtain
$$\partial_tW+\div \Big(-[\sum_l\sum_id^i_{kl}\partial_lc_i]_k+\sum_iu_ic_i\Big)=\sum_if_i\leq W\,b_0,$$
hence
$$\partial_tW-\sum_k\partial_k\Big[\sum_l\sum_id^i_{kl}\partial_lc_i\Big]+\div \Big(\sum_iu_ic_i\Big)\leq W\,b_0.$$
Equivalently,
\begin{equation}\label{Sigmas}
\partial_tW-\sum_{k}\partial_{k}\Big[\sum_l\sum_i \Big( \partial_l (d^i_{kl}c_i)-(\partial_ld^i_{kl})c_i\Big) \Big]+\div (\sum_iu_ic_i)\leq W\,b_0.
\end{equation}
We set
$$A_{kl}=\frac{1+\sum_id^i_{kl}c_i}{1+\sum_jc_j},\;\;B=[B_k]_k,\;\; B_k=\sum_i\frac{c_i}{1+\sum_jc_j}\sum_l\partial_l(d^i_{kl}),\;\; U=\sum_i\frac{c_i}{1+\sum_j c_j}u_i.$$
Then (\ref{Sigmas}) may be rewritten as
\begin{equation}\label{Waniso}
\partial_t W-\sum_{k,l}\partial_{kl}(A_{kl}W)+\div  (B\,W)+\div (U\,W)\leq W\,b_0.
\end{equation}
The sum of the boundary conditions leads to
\begin{equation}\label{Wbdyaniso}
 -\sum_{k,l}\partial_l(A_{kl}W)\nu_k+W(B+U)\cdot\nu=0.
 \end{equation}

\noindent
Note the following estimates:
\begin{equation}\label{sumestimates}
\left.
\begin{array}{l}
\|U\|_{L^\infty(0,T;L^r(\Omega))}\leq \max_i\|\frac{c_i}{W}\|_{L^\infty(Q_T)}\sum_i\|u_i\|_{L^\infty(0,T;L^r(\Omega))}\leq \sum_i\|u_i\|_{L^\infty(0,T;L^r(\Omega))},\\[0.5ex]
\|B\|_{L^\infty(0,T;L^r(\Omega))}\leq \max_k \sum_{i,l}\|\partial_l(d^i_{kl})\|_{L^\infty(0,T;L^r(\Omega))},\\[0.5ex]
\forall \xi,\; \sum_{k,l}A_{kl}\xi_k\xi_l=\frac{|\xi|^2+\sum_i\sum_{k,l}d^i_{kl}\xi_k\xi_lc_i}{1+\sum_jc_j}\geq \frac{|\xi|^2+\sum_i\alpha_i|\xi|^2c_i}{1+\sum_jc_j}\geq \min\{1,\min_i\alpha_i\}|\xi|^2,\\[0.5ex]
A_{kl}\in C(\overline{Q_T}),\; \|A_{kl}\|_{L^\infty(Q_T)}\leq \max\{1, \max_i\|d"_{kl}"\|_{L^\infty(Q_T)}\}.
\end{array}
\right\}
\end{equation}
Let us introduce the dual problem of (\ref{Waniso}-\ref{Wbdyaniso}), namely
\begin{equation}\label{dualaniso}
\left.
\begin{array}{l}
-\big[\partial_t\Phi+\sum_{k,l}A_{kl}\partial_{kl}\Phi+(B+U)\nabla\Phi+b_0\Phi\big]=\Theta,\\[0.5ex]
\nabla\Phi\cdot[\sum_lA_{kl}\nu_l]_k=0  \mbox{ on } \Sigma_T,\\[0.5ex]
\Phi(T)=0.
\end{array}
\right\}
\end{equation}
Then, we use the following Alexandrov-Bakelman-Pucci-Krylov-Safonov type of estimate.
\begin{lemma} [\cite{NU}]\label{krysaf}
There exist $K>0$ depending only on the data in (\ref{sumestimates}) and on $T, b_0$ such that
$$\|\Phi\|_{L^\infty(Q_T)}\leq K\|\Theta\|_{L^{N+1}(Q_T)},$$
\end{lemma}
\noindent
where $\Phi$ is the weak solution (in the sense given above) of \eqref{dualaniso}.\\[1ex]
\noindent
Assume $\Theta\geq 0$ in (\ref{dualaniso}) so that $\Phi\geq 0$. Multiplying (\ref{Waniso}) by $\Phi$ and integrating by parts lead to
$$-\int_\Omega W(0)\Phi(0)+\int_{\Sigma_T}\sum_{k,l}A^{kl}W(\partial_k \Phi)\nu_l-\partial_l(A^{kl}W)\nu_k\Phi+W(B+U)\Phi\cdot\nu+\int_{Q_T}\Theta\,W\leq 0.$$
Using the boundary conditions on $W$ and $\Phi$, this gives
$$\int_{Q_T}\Theta\,W\leq \int_\Omega W(0)\Phi(0).$$
Since Lemma~\ref{krysaf} implies $\|\Phi(0)\|_{L^\infty(\Omega)}\leq K\|\Theta\|_{L^{N+1}(Q_T)}$ for the continuous
$\Phi$, we deduce
$$\|W\|_{L^{(N+1)/N}(Q_T)}\leq K\|W(0)\|_{L^1(\Omega)}$$
by duality. The estimate of Proposition \ref{LN+1} follows.
\cqfd
\begin{remark}{\rm Actually, according to the results in \cite{N}, even the $C^\alpha(Q_T)$-norm of $\Phi$ can be estimated in terms of the $L^{N+1}(Q_T)$-norm of $\Theta$ in Problem (\ref{dualaniso}), at least if $r$ in \eqref{sumestimates} satisfies $r\geq N+2$. Therefore, the mapping $\Theta\in L^{N+1}(Q_T)\to \Phi\in L^\infty(Q_T)$ is then even compact. By duality, one can then deduce compactness properties of $W$ in $L^{(N+1)/N}(Q_T)$.
}
\end{remark}
\section*{Acknowledgement}
The second author was partly supported by the Center of Smart Interfaces, TU Darmstadt and Rennes Metropole.
\bibliography{literature}{}

\begin{thebibliography}{10}

\bibitem{amann89}
H.~Amann.
\newblock Dynamic theory of quasilinear parabolic systems. {III}. {G}lobal
  existence.
\newblock {\em Math. Z.}, 202(2):219--250, 1989.

\bibitem{amann93}
H.~Amann.
\newblock Nonhomogeneous linear and quasilinear elliptic and parabolic boundary
  value problems.
\newblock In {\em Function spaces, differential operators and nonlinear
  analysis ({F}riedrichroda, 1992)}, volume 133 of {\em Teubner-Texte Math.},
  pages 9--126. Teubner, Stuttgart, 1993.

\bibitem{Bothe-MS}
D.~Bothe.
\newblock On the {M}axwell-{S}tefan approach to multicomponent diffusion.
\newblock In {\em Progress in Nonlinear Differential Equations and their
  Applications}, pages 81--93. Springer, Basel, 2011.

\bibitem{BotheDreyer}
D.~Bothe and W.~Dreyer.
\newblock Continuum thermodynamics of chemically reacting fluid mixtures.
\newblock {\em Acta Mechanica}, pages 1757--1805, 2015.

\bibitem{BFPR12b}
D.~Bothe, A.~Fischer, M.~Pierre, and G.~Rolland.
\newblock Global existence for diffusion-electromigration systems in space
  dimension three and higher.
\newblock {\em Nonlinear Analysis A: Theory Methods and Applications},
  99:152--166, 2014.

\bibitem{BPR12}
D.~Bothe, M.~Pierre, and G.~Rolland.
\newblock Cross-diffusion limit for a reaction-diffusion system with fast
  reversible reaction.
\newblock {\em Comm. in P.D.E}, 37:1940--1966, 2012.

\bibitem{BotheRolland14}
D.~Bothe and G.~Rolland.
\newblock Global existence for a class of reaction-diffusion systems with mass
  action kinetics and concentration-dependent diffusivities.
\newblock {\em Acta Appl. Math.}, 139(1):25--57, 2015.

\bibitem{Brezis}
Ha{\"{\i}}m Brezis.
\newblock {\em Functional analysis, {S}obolev Spaces and Partial Differential
  Equations}.
\newblock Universitext. Springer, New York, 2011.

\bibitem{CantrellAdvection}
R.S. Cantrel, C.~Cosner, and Y.~Lou.
\newblock Advection-mediated coexistence of competing species.
\newblock {\em Proc. Royal Soc. Edinburgh}, 137:497--512, 2007.

\bibitem{CantrellCosner}
R.S. Cantrell and C.~Cosner.
\newblock {\em Spatial Ecology via Reaction-Diffusion Equations}.
\newblock John Wiley \& Sons, 2003.

\bibitem{Levine}
T.~Chen, H.A. Levine, and Sacks P.E.
\newblock Analysis of a convective reaction-diffusion equation.
\newblock {\em Nonlinear Anal., Theor., Meth. and Appl.}, 12:1349--1370, 1988.

\bibitem{dhp}
R.~Denk, M.~Hieber, and J.~Pr{\"u}ss.
\newblock Optimal {$L^p-L^q$}-regularity for parabolic problems with
  inhomogeneous boundary data.
\newblock {\em Math. Z.}, 257:193--224, 2007.

\bibitem{Evans}
L.~C. Evans.
\newblock {\em Partial differential equations}.
\newblock Providence, RI: American Mathematical Society, 1998.

\bibitem{FeireislTrivisa}
E.~Feireisl, H.~Petzeltova, and K.~Trivisa.
\newblock Multicomponent reactive flows: Global-in-time existence for large
  data.
\newblock {\em Comm. Pure Appl. Anal.}, 7:1017--1047, 2008.

\bibitem{fischer13}
A.~Fischer.
\newblock {\em Well-posedness and asymptotic behaviour in reactive and
  electro-kinetic flow processes}.
\newblock {PhD} thesis, Technische Universit{\"a}t Darmstadt, 2013.

\bibitem{GG96}
H.~Gajewski and K.~Gr{\"o}ger.
\newblock Reaction-diffusion processes of electrically charged species.
\newblock {\em Mathematische Nachrichten}, 177:109--130, 1996.

\bibitem{GoudonVasseur}
T.~Goudon and A.~Vasseur.
\newblock Regularity analysis for systems of reaction-diffusion equations.
\newblock {\em Ann. Sci. Ec. Norm. Sup.}, 43:117--142, 2010.

\bibitem{hollis}
S.~L. Hollis, R.~H. Martin, Jr., and M.~Pierre.
\newblock Global existence and boundedness in reaction-diffusion systems.
\newblock {\em SIAM J. Math. Anal.}, 18(3):744--761, 1987.

\bibitem{HollisMorgan}
S.~L. Hollis and J.~Morgan.
\newblock Global existence and asymptotic decay for systems of convective
  reaction-diffusion equations.
\newblock {\em Nonlinear Anal., Theor., Meth. and Appl.}, 17:725--739, 1991.

\bibitem{Kraeutle}
K.~Kr{\"a}utle.
\newblock Existence of global solutions of multicomponent reactive transport
  problems with mass action kinetics in porous media.
\newblock {\em J. Appl. Ana. Comp.}, 1:497--515, 2011.

\bibitem{LSU}
O.~A. Lady{\v{z}}enskaja, V.~A. Solonnikov, and N.~N. Ural{'}ceva.
\newblock {\em Linear and quasilinear equations of parabolic type}.
\newblock Translated from the Russian by S. Smith. Translations of Mathematical
  Monographs, Vol. 23. American Mathematical Society, Providence, R.I., 1967.

\bibitem{Mo}
J.~Morgan.
\newblock Global existence for semilinear parabolic systems.
\newblock {\em SIAM J. Math. Anal.}, 20(5):1128--1144, 1989.

\bibitem{N}
A.I. Nazarov.
\newblock H{\"o}lder estimates for bounded solutions of problems with an
  oblique derivative for parabolic equations with nondivergent structure.
\newblock {\em J. Sov. Math. translated from Probl. Mat. Anal., 11 (1990),
  37-46 (Russian)}, 64(6):1247--1252, 1993.

\bibitem{NU}
A.I. Nazarov and N.N. Uraltseva.
\newblock The oblique boundary-value problem for a quasilinear parabolic
  equation.
\newblock {\em J. Math. Sc. translated from ZNS POMI, 200 (1992), 118-131
  (Russian)}, 77(3):3212--3220, 1995.

\bibitem{pierre10}
M.~Pierre.
\newblock Global existence in reaction-diffusion systems with control of mass:
  a survey.
\newblock {\em Milan J. Math.}, 78(2):417--455, 2010.

\bibitem{PSch}
M.~Pierre and D.~Schmitt.
\newblock Blowup in reaction-diffusion systems with dissipation of mass.
\newblock {\em SIAM Rev.}, 42(1):93--106 (electronic), 2000.

\bibitem{pruss}
J.~Pr{\"u}ss.
\newblock Maximal regularity for evolution equations in {$L_p$}-spaces.
\newblock {\em Conf. Semin. Mat. Univ. Bari}, 285:1--39 (2003), 2002.

\bibitem{rolland12}
G.~Rolland.
\newblock {\em Global existence and fast-reaction limit in reaction-diffusion
  systems with cross effects}.
\newblock {PhD} thesis, ENS Cachan-Bretagne, 2012.

\bibitem{simon86}
J.~Simon.
\newblock Compact sets in the space {$L^p(0,T;B)$}.
\newblock {\em Annali di Matematica Pura ed Applicata}, 146:65--96, 1987.

\bibitem{Texier}
R~Texier-Picard and V.~Volpert.
\newblock Reaction-diffusion-convection problems in unbounded cylinders.
\newblock {\em Revista Matem\`atica Complutense}, 16:223--276, 2003.

\end{thebibliography}
\bibliographystyle{plain}
\end{document}